\documentclass[11pt,a4paper]{article}
\usepackage{amssymb,amsmath,amsfonts}
\usepackage{url}
\usepackage{epsfig}
\usepackage{hyperref}
\usepackage{srcltx}
\usepackage{xcolor}
\usepackage{lineno}

\def\mK{{\mathfrak{T}}}

\newcommand{\qed}{\hfill $\Box$ \nr \medskip}
\def\lin{{\mathrm{lin}}}
\def\poly{{\mathrm{poly}}}

\def\three?{3}
\def\four?{4}
\def\ten?{10}

\def\bK{{\mathbf{K}}}

\def\Argmin{\mathop{\hbox{\rm Argmin}}}
\def\beq{\begin{equation}}
\def\eeq{\end{equation}}

\def\norm2to2{{\|\cdot\|_{2,2}}}
\def\Prob{\hbox{\rm Prob}}

\def\bE{{\mathbf{E}}}

\def\inter{\hbox{\rm  int}}

\def\S{{\mathbf{S}}}
\def\Diag{\hbox{\rm  Diag}}
\def\Prob{\hbox{\rm  Prob}}

\def\Opt{\hbox{\rm Opt}}

\def\Card{{\mathop{\hbox{\rm  Card}}}}

\def\n{{\hbox{\rm  n}}}

\def\Tr{{\mathop{\hbox{\rm  Tr}}}}

\def\cB{{\cal B}}

\def\cH{{\cal H}}

\def\cM{{\cal M}}
\def\cN{{\cal N}}

\def\cR{{\cal R}}
\def\cS{{\cal S}}
\def\cT{{\cal T}}

\def\cW{{\cal W}}
\def\cX{{\cal X}}

\def\cZ{{\cal Z}}

\def\S{{\cal S}}

\def\Argmin{\mathop{\hbox{\rm  Argmin}}}

\def\S{{\mathbf{S}}}

\def\bK{{\mathbf{K}}}

\def\bS{{\mathbf{S}}}

\def\bQ{{\mathbf{Q}}}

\def\qed{\ \hfill$\square$\par\smallskip}

\def\mypict3{\epsfxsize=220pt\epsfysize=80pt\epsfbox}

\def\bR{{\mathbf{R}}}

\def\argmin{\mathop{\hbox{\rm argmin}}}
\def\cH{{\cal H}}
\def\Col{{\hbox{\rm Col}}}

\def\Risk{{\hbox{\rm Risk}}}
\long\def\comment#1{}

\def\poly{{\mathrm{poly}}}

\newtheorem{proposition}{Proposition}[section]

\newtheorem{theorem}{Theorem}[section]
\newtheorem{example}{Example}[section]

\def\lin{{\hbox{\rm\tiny lin}}}

\def\three?{3}
\def\four?{4}
\def\ten?{10}

\def\bK{{\mathbf{K}}}

\def\Argmin{\mathop{\hbox{\rm Argmin}}}
\def\beq{\begin{equation}}
\def\eeq{\end{equation}}

\def\norm2to2{{\|\cdot\|_{2,2}}}
\def\Prob{\hbox{\rm Prob}}

\def\bE{{\mathbf{E}}}

\def\inter{\hbox{\rm  int}}

\def\S{{\mathbf{S}}}
\def\Diag{\hbox{\rm  Diag}}
\def\Prob{\hbox{\rm  Prob}}

\def\Opt{\hbox{\rm Opt}}

\def\Card{{\mathop{\hbox{\rm  Card}}}}

\def\n{{\hbox{\rm  n}}}

\def\Tr{{\mathop{\hbox{\rm  Tr}}}}

\def\cB{{\cal B}}

\def\cH{{\cal H}}

\def\cM{{\cal M}}
\def\cN{{\cal N}}

\def\cR{{\cal R}}
\def\cS{{\cal S}}
\def\cT{{\cal T}}

\def\cW{{\cal W}}
\def\cX{{\cal X}}

\def\cZ{{\cal Z}}

\def\S{{\cal S}}

\def\Argmin{\mathop{\hbox{\rm  Argmin}}}

\def\S{{\mathbf{S}}}

\def\bK{{\mathbf{K}}}

\def\bS{{\mathbf{S}}}

\def\bQ{{\mathbf{Q}}}

\def\qed{\ \hfill$\square$\par\smallskip}

\def\mypict3{\epsfxsize=220pt\epsfysize=80pt\epsfbox}

\def\bR{{\mathbf{R}}}

\def\argmin{\mathop{\hbox{\rm argmin}}}
\def\cH{{\cal H}}
\def\Col{{\hbox{\rm Col}}}

\def\Risk{{\hbox{\rm Risk}}}
\newcommand{\hide}[1]{{}}
\newcommand{\aic}[2]{{#2}}

\newcommand{\be}{\begin{eqnarray}}
\newcommand{\ee}[1]{\label{#1}\end{eqnarray}}
\newcommand{\nn}{\nonumber \\}
\newcommand{\ese}{\end{eqnarray*}}
\newcommand{\bse}{\begin{eqnarray*}}
\newcommand{\rf}[1]{~(\ref{#1})}
\def\SG{{\cal SG}}

\def\half{\tfrac{1}{2}}
\def\four{\tfrac{1}{4}}
\newcommand{\wh}[1]{{\widehat{#1}}}

\def\wt#1{\widetilde{#1}}
\def\ov#1{\overline{#1}}

\begin{document}
\title{On robust recovery of signals from indirect observations}
\author{Yannis Bekri\thanks{LJK, Universit\'e Grenoble Alpes, 700 Avenue Centrale,  38401 Domaine Universitaire de Saint-Martin-d'Hères, France {\tt yannisbekri@hotmail.com}} \and
Anatoli Juditsky
\thanks{LJK, Universit\'e Grenoble Alpes, 700 Avenue Centrale,  38401 Domaine Universitaire de Saint-Martin-d'Hères, France {\tt anatoli.juditsky@univ-grenoble-alpes.fr}}
\and Arkadi Nemirovski
\thanks{Georgia Institute
 of Technology, Atlanta, Georgia
30332, USA, {\tt nemirovs@isye.gatech.edu}\newline
This work was supported by Multidisciplinary Institute in Artificial intelligence MIAI {@} Grenoble Alpes (ANR-19-P3IA-0003).}}
\date{}
\maketitle
\begin{abstract}
We consider an {\em uncertain linear inverse problem} as follows. Given observation
$\omega=Ax_*+\zeta$ where $A\in \bR^{m\times p}$ and $\zeta\in \bR^{m}$ is observation noise, we want to recover unknown signal $x_*$, known to belong to a convex set ${\cal X}\subset \bR^{n}$. As opposed to the ``standard'' setting of such problem, we suppose that the model noise $\zeta$ is ``corrupted''---contains an uncertain (deterministic dense or singular) component. Specifically, we assume that $\zeta$ decomposes into $\zeta=N\nu_*+\xi$ where $\xi$ is the random noise and $N\nu_*$ is the ``adversarial contamination'' with known  $\cN\subset \bR^n$ such that $\nu_*\in \cN$ and $N\in \bR^{m\times n}$.
We consider two ``uncertainty setups'' in which $\cN$ is either a convex bounded set or is the set of sparse vectors (with at most $s$ nonvanishing entries).  We analyse the performance of  ``uncertainty-immunized'' {\em polyhedral estimates}---a particular class of nonlinear estimates as introduced in \cite{juditsky2020polyhedral,PUP}---and show how ``presumably good'' estimates of the sort may be constructed in the situation where the signal set is an {\em ellitope} (essentially, a symmetric convex set delimited by quadratic surfaces) by means of efficient convex optimization routines.

\end{abstract}
\section{Situation  and goals} \label{secsitgoals}
\subsection{Introduction}
Since the term was coined in the 1950's, the problem of robust estimation has received much attention in the classical statistical literature.\footnote{It is impossible to give an overview of the existing literature on robust estimation, and we do not try to do it here; for the ``classical'' framework one may refer to early references in \cite{tukey1960survey}, the foundational manuscript \cite{Huber1981}, or a recent survey \cite{yu2017robust}.} In this paper our focus is on robust estimation of a signal from indirect linear observations.
Specifically, suppose that our objective is to recover a linear image $w_*=Bx_*$ of unknown signal $x_*$, known to belong to a given convex set $\cX\subset \bR^p$, given $B\in \bR^{q\times p}$, $A\in \bR^{m\times p}$, and a noisy observation
\be\omega=Ax_*+\eta_*+\xi \in \bR^m
\ee{modintro}
of $x_*$,  perturbed by a mixed---random-and-deterministic noise $\xi+\eta_*$. Here $\xi$ is the random noise component, while  $\eta_*$ is the ``adversarial'' deterministic noise.
Recently, this problem attracted much attention in the context of robust recovery of sparse ({with} at most $s\ll p$ nonvanishing entries) signal $x_*$.
In particular, robust sparse regression with
an emphasis on contaminated design was investigated in \cite{chen2013robust,balakrishnan2017computationally,diakonikolas2019efficient,liu2020high,dalalyan2019outlier,minsker2024robust}; methods based on penalizing the
vector of outliers were studied {in} \cite{foygel2014corrupted}, see also \cite{bruce1994denoising,sardy2001robust}.

In this paper, our {emphasis} is on {specific} assumptions about the structure of the signal $x_*$ to {be recovered}  and {on the} contamination $\eta_*$.
Namely, we assume that the set $\cX$ of signals is an {\em ellitope}---a convex compact symmetric w.r.t. the origin  subset of $\bR^p$ delimited by quadratic surfaces.\footnote{See \cite[Section 4.2.1]{PUP} or Section \ref{secellit} below; as of now, an instructive example of ellitope is an  intersection of finite family of ellipsoids/elliptic cylinders with common center.}
Our interest for ellitopes is motivated by {the fact that these signal sets are well suited for the} problem of estimating unknown signal $x_*$ from observation \rf{modintro} in the {Gaussian} no-nuisance case  ($\eta_*=0$,
$\xi\sim\cN(0,\sigma^2 I_m)$). Specifically, let us consider  linear estimate $\wh w_{\lin}(\omega)=G^T_{\lin}\omega$ and {\em polyhedral estimate} $\wh w_\poly(\omega)=B\wh x_\poly(\omega)$ where
\[
\wh x_\poly(\omega)\in \Argmin_{x\in \cX} \|G^T_{\poly}(Ax-\omega)\|_\infty
\]
of $w_*$. Let $\cX$ be an ellitope, and let the estimation error be measured in a co-ellitopic norm $\|\cdot\|$ (i.e., such that the unit ball $\cB_*$ of the norm $\|\cdot\|_*$ conjugate to $\|\cdot\|$ is an ellitope). In this situation, one can point out (cf. \cite{JudNem2018,juditsky2020polyhedral,PUP})  efficiently computable {\em contrast matrices} $ G_{\lin}\in \bR^{m\times q}$ and $G_{\poly}\in\bR^{m\times m}$ such that estimates $\wh w_{\lin}(\cdot)$ and $\wh w_\poly(\cdot)$  attain nearly minimax-optimal performance.

We suppose that the adversarial perturbation $\eta_*$ has  a special structure: we are given a ``nuisance set'' $\cN\subset\bR^n$ such that $\nu_*\in \cN$ and a $m\times n$ matrix $N$ such that $\eta_*=N\nu_*$. We consider two types of assumptions about $\cN$:  $\cN$ is either 1) a (nonempty) compact convex set, or, more conventionally, 2) $\cN$ is the set of sparse disturbances (with at most $s\leq n$ nonvanishing components). Our focus is on the design and performance analysis of the ``uncertainty immunized'' polyhedral estimate $\wh w(\omega)$ of $w_*=Bx_*$  in the presence of the contaminating signal, and solving the problem in the first case leads to a ``presumably good'' solution for the second.

We would like to emphasize the principal feature of the approach we promote: in this paper, $A,\,B$, and $N$ are ``general'' matrices of appropriate dimensions, while $\cX$ and $\cN$  are rather general sets. As a consequence, we adopt here an ``operational'' approach\footnote{As opposed to the classical ``descriptive'' approach to solving the estimation problem in question via deriving and optimizing \aic{}{w.r.t. estimate parameters} closed form analytical expressions for the risk of a candidate estimate.}
initiated in \cite{Don95} and further developed in \cite{jn09,l2estimation,juditsky2020polyhedral,PUP},
within which both the estimates and their risks are yielded by efficient computation, rather than by an explicit analytical analysis, seemingly impossible under the circumstances. The term ``efficient'' in the above is essential, and is also responsible for principal limitations of the results to follow. First of all, it imposes restrictions on the structure of the set of signals of interest and on the norm {quantifying} the estimation error. As it is shown in \cite{PUP}, the maximum of a quadratic form over an ellitope admits a ``reasonably tight'' efficiently computable upper bound, leading to tight bounds {on the risk of linear and polyhedral estimates when the signal set is an ellitope}. Furthermore, while in the case of convex compact set $\cN$ of contaminations, constructing risk bounds for the polyhedral estimate $\wh w_G(\omega)$ associated with a {\em given} contrast matrix $G$ is possible under rather weak assumptions about the nuisance set $\cN$ (essentially, the computational tractability of this set is sufficient), the fundamental problem of {\em contrast synthesis}---minimizing these bounds over contrast matrices---allows for efficiently computable solution only when $\cN$ is either {an} ellitope itself, or is a ``co-ellitope'' (the polar of an ellitope).

To {complete} this section, we would like to mention another line of research on the problem of estimating signal $x_*$ from observation \rf{modintro} under purely deterministic disturbance (case of $\xi=0$), the standard problem of optimal recovery \cite{MicRiv77,MicRiv84} and guaranteed estimation under uncertain-but-bounded perturbation \cite{chernousko1993state,fogel1982value,granichin2003randomized,kurzhanski1989identification,kurzhanski1997ellipsoidal,milanese1991optimal,schweppe1973uncertain}.
The present work may be seen as an attempt to extend the corresponding framework to the case in which both deterministic and random observation noises are present.

%

\paragraph{Organization of the {paper}.}
We introduce the exact statement of the {the} estimation problem {to be considered} and the entities which are relevant for the analysis to follow in Section \ref{theproblem}.
Analysis and design of the polyhedral estimate in the case of uncertain-but-bounded contamination are presented in Section \ref{sec:bond}. Then in Section \ref{13estimate} we describe the application of the proposed framework to the case of (unbounded) singular contamination using the sparse {model} of the nuisance vector.
\aic{}{Finally, we recall some results on $\ell_1$ recovery used in the paper in appendix \ref{sectspars}.}


\subsection{The problem}\label{theproblem}
The estimation problem we are interested in is as follows:
\begin{quote}
Recall that we are given observation (cf. \rf{modintro})
\be
\omega=Ax_*+N \nu_*+\xi \in \bR^m
\ee{omega}
where
\begin{itemize}
\item $N\in\bR^{m\times n}$, $A\in\bR^{m\times p}$ are given matrices,
\item $\nu_*\in\bR^n$ is unknown {\em nuisance signal}, $\nu_*\in \cN$, a known subset of $\bR^n$,
\item $x_*$ is unknown signal of interest known to belong to a given convex compact set $\cX\subset\bR^p$ symmetric w.r.t. the origin,
\item $\xi\sim P_\xi$ is a random observation noise.
\end{itemize}
Given $\omega$, {\em our objective} is to recover the linear image $w_*=Bx_*$ of $x_*$, $B$ being a given ${q\times p}$ matrix.
\end{quote}
Given $\epsilon\in(0,1)$, we quantify the quality of the recovery $\wh w(\cdot)$ by its $\epsilon$-risk\footnote{The $\epsilon$-risk of an estimate depends, aside of $\epsilon$ and the estimate, on the ``parameters'' $\|\cdot\|$, $\cX$, $\cN$; these entities will always be specified by the context, allowing us to omit mentioning them in the notation $\Risk_\epsilon[\cdot]$.}
$$
\Risk_\epsilon[\widehat{w}]=\inf\left\{\rho:
\,\Prob_{\xi\sim P_\xi}\{\xi:\,\|Bx_*-\widehat{w}(\omega)\|
>\rho\}\leq \epsilon\quad
\forall (\hbox{$\nu_*\in \cN$ and $x_*\in\cX$})\right\}
$$
where $\|\cdot\|$ is a given norm.

\paragraph{Observation noise assumption.}
In the sequel, we assume that the observation noise $\xi$ is sub-Gaussian with parameters $(0,\sigma^2I_m)$, that is,
\[
\bE\left\{e^{h^T\xi}\right\}\leq \exp\left(\half \sigma^2\|h\|_2^2\right).
\]
\subsection{Ellitopes}\label{secellit}
Risk analysis of a candidate polyhedral estimate heavily depends on the geometries of the signal set $\cX$ and norm $\|\cdot\|$. In the sequel, we restrict
ourselves to the case where $\cX$ and the {\sl polar} $\cB_*$ of the unit ball of $\|\cdot\|$ are {\sl ellitopes}.
\par 
By definition \cite{JudNem2018,PUP}, {\em a basic ellitope} in $\bR^n$ is a set of the form
\be   \cX=\{x\in\bR^n:\,\exists z\in\bR^k,\,t\in\cT:\,z^TT_\ell z\leq t_\ell,\,\ell\leq L\},
\ee{ell1} where $T_\ell\in\bS^k$, $T_\ell\succeq0$, $\sum_\ell T_\ell \succ0$, and $\cT\subset\bR^L_+$ is a convex compact set with a nonempty interior which is monotone: whenever $0\leq t'\leq t\in\cT$ one has $t'\in\cT$. An ellitope is an image of a basic ellitope under a linear mapping.
We refer to $L$ as {\em ellitopic dimension} of $\cX$.
\par
Clearly, every ellitope is a convex compact set symmetric w.r.t. the origin;  a {\em basic ellitope}, in addition, has a nonempty interior.
\paragraph{Examples.}~\\ {\bf A.} Bounded intersection $\cX$ of $L$ centered at the origin ellipsoids/elli\-p\-tic cylinders $\{x\in \bR^n:\,x^T T_\ell x\leq1\}$ [$T_\ell\succeq0$] is a basic ellitope:
\[
\cX=\{x\in \bR^n:\exists t\in\cT:=[0,1]^L: x^T T_\ell x\leq t_\ell,\,\ell\leq L\}
\]
In particular, the unit box $\{x\in \bR^n:\|x\|_\infty\leq1\}$ is a basic ellitope.\\
{\bf B.} A $\|\cdot\|_p$-ball in $\bR^n$ with $p\in[2,\infty]$ is a basic ellitope:
\[
\{x\in\bR^n:\|x\|_p\leq1\} =\big\{x:\exists t\in\cT=\{t\in\bR^n_+,\|t\|_{p/2}\leq 1\}:\underbrace{ x_\ell^2}_{x^T T_\ell x}\leq t_\ell,\,\ell\leq n\big\}.
\]
Ellitopes admit fully algorithmic ``calculus''---this family is closed with respect to basic operations preserving convexity
and symmetry w.r.t. the origin, e.g., taking finite intersections, linear images, inverse images under linear embedding, direct products, arithmetic summation
 (for details, see \cite[Section 4.6]{PUP}).
\paragraph{Main assumption.}
We assume from now on that the signal set $\cX$ and the {\sl polar} $\cB_*$ of the unit ball of the norm $\|\cdot\|$ are {\em basic ellitopes:}\footnote{The results to follow straightforwardly extend to the case where $\cX$ and $\cB_*$ are ``general'' ellitopes; we assume them to be basic to save notation.}
\begin{subequations}\label{cXcB}
\begin{align}\label{elliX}
\cX&=\{x\in\bR^n:\exists t\in \cT: x^TT_kx\leq t_k,k\leq K\},\\
\cB_*&=\{y\in\bR^q:\exists s\in\cS: y^TS_\ell y\leq s_\ell,\ell\leq L\}\label{elliB}
\end{align}
\end{subequations}
where $\cT\subset\bR^K_+$, $\cS\in\bR^L_+$ are monotone convex compact sets with nonempty interiors, $T_k\succeq0,\sum_kT_k\succ0$, $S_\ell\succeq0$, and $\sum_\ell\S_\ell\succ0$.\par
\paragraph{Notation.}
In the sequel, \aic{}{$[A;B]$ and $[A,B]$ stand for vertical and horizontal concatenation of matrices $A$ and $B$ of appropriate dimensions;
notation  $A\succeq B$ ($A\succ B$) means that matrix $A-B$ is positive semidefinite (respectively, positive definite).}
We denote $\n[G]$ the maximum of
Euclidean norms of the columns of a matrix $G$.
\section{Uncertain-but-bounded nuisance}\label{sec:bond}
In this  section, we consider the case of {\sl uncertain-but-bounded} nuisance. Specifically, we assume that $\cN\subset\bR^n$ is a convex compact set, symmetric w.r.t. the origin, and specify $\pi(\cdot)$ as the semi-norm on $\bR^m$ given by
\[\pi(h)=\sup_{u}\left\{(Nu)^Th:\,u\in \cN\right\}.
\]
\subsection{Bounding the $\epsilon$-risk of polyhedral estimate}
In this section, a polyhedral estimate is specified by $m\times I$ {\sl contrast matrix} $G$ and is as follows: given observation $\omega$ (see (\ref{omega})), we find an optimal solution $\wh x_G(\omega)$ to the
(clearly solvable)  optimization problem
\be
\min_{x,\nu}\big\{\|G^T(Ax+N\nu-\omega)\|_\infty:\,x\in \cX,\,\nu\in \cN\big\}.
\ee{1stG}
Given a $m\times I$ matrix $G=[g_1,..,g_I]$, let  $\Xi_\epsilon[G]$ be the set of all realizations of $\xi$ such that
\be|[g_i^T\xi]_i|\leq\underbrace{\sigma\sqrt{2\ln \left[{2I}/\epsilon\right]}}_{=:\varkappa(\epsilon)}\|{g_i}\|_2,\quad \forall i\leq I.
\ee{Xi_e}
Note that
\begin{equation}\label{neweq2}
\Prob_{\xi\sim\SG(0,\sigma^2 I_m)}\left\{\xi\not\in\Xi_{\epsilon}(G)\right\}\leq \epsilon.
\end{equation}
Indeed, we have $\bE_\xi\left\{e^{\gamma g^T\xi}\right\}\leq e^{{1\over 2}\gamma^2\|g\|_2^2\sigma^2}$, implying that for all $a\geq0$,
\[
\Prob\{g^T\xi>a\}\leq\inf_{\gamma>0}\exp\left\{\half \gamma^2\|g\|_2^2\sigma^2-\gamma a\right\}=\exp\left\{-\half{a^2\sigma^2\|g\|_2^2}\right\},
\]
so that
\[\Prob\left\{\exists i\leq I: |g_i^T\xi|>\varkappa(\epsilon)\|g_i\|_2\right\} \leq 2I\exp\left\{-\tfrac{\varkappa^2(\epsilon)}{2\sigma^2}\right\}\leq\epsilon.
\]

\begin{proposition}\label{pro:G2} Given a $m\times I$ contrast matrix $G=[g_1,...,g_I]$,  consider the optimization problem
\begin{align}\label{f_G2}
\Opt[G]&=\min\limits_{\lambda,\mu,\gamma}\bigg\{
\overbrace{\phi_\cS(\lambda)+4\phi_\cT(\mu)+4\psi^2[G]\sum_i\gamma_i}^{=:f_G(\lambda,\mu,\gamma)}:\;
\lambda\geq0,\;\mu\geq0,\,\gamma\geq0,\\
&\quad\quad\quad\quad
\left.\left[\begin{array}{c|c}\sum_\ell{\lambda}_\ell S_\ell&{1\over 2}B\\\hline
{1\over 2}B^T&
\sum_k \mu_kT_k+
A^T\left[\sum_i\gamma_ig_ig_i^T\right]A
\cr\end{array}\right]\succeq0\right\}\nonumber
\end{align}
where
$$
\psi[G]=\max_i\pi(g_i)+\varkappa(\epsilon)\n[G]
$$
and, from now on, for a nonempty compact set $\cZ\subset \bR^N$
\[
\phi_{\cZ}(\zeta):=\max_{z\in\cZ}\zeta^T t
\] is the support function of $\cZ$.
Let $(\lambda,\mu,\gamma)$ be a feasible solution to the problem in \rf{f_G2}. Then
\[\Risk_\epsilon[\wh w_G]\leq f_G(\lambda,\mu,\gamma),
\]
i.e., the $\epsilon$-risk of the estimate $\wh w_G$ is upper bounded with $f_G(\lambda,\mu,\gamma)$.
\end{proposition}
{\bf Proof.} Let us fix $\xi\in\Xi_\epsilon[G]$, $x_*\in\cX$, and $\eta_*\in\cN$. 
Let also $\wh x=\wh x_G(\omega)$ be the $x$-component of some optimal solution $[\wh x;\wh\nu]$, $\wh\nu\in\cN$,  to (\ref{1stG}) and, finally, let $\Delta=\wh x-x_*$. Observe that  $[x,\nu]=[x_*,\nu_*]$ is feasible for \rf{1stG} and $\|G^T[Ax_*+N\nu_*-\omega]\|_\infty=\|G^T\xi\|_\infty\leq \varkappa(\epsilon)\n[G]$, implying that $\|G^T[A\wh x+N\wh \nu-\omega]\|_\infty \leq\varkappa(\epsilon)\n[G]$ as well. Therefore,
$$
\|G^TA\Delta\|_\infty\leq 2\varkappa(\epsilon)\n[G]+\|G^TN[\wh\nu-\nu_*]\|_\infty.
$$
Taking into account that $\wh\nu,\nu_*\in\cN$, we have $\|G^TN[\wh\nu-\nu_*]\|_\infty\leq2\max_i\pi(g_i)$, and we arrive at
\begin{equation}\label{role}
|g_i^TA\Delta|\leq 2\psi[G],\,i=1,...,I.
\end{equation}
Now,  we have $\Delta\in 2\cX$, that is, for some $t\in \cT$ and all $k$ it holds $\Delta^TT_k\Delta\leq 4t_k$, and let $v\in\cB_*$, so that for some $s\in\cS$ for all $\ell$ it holds $v^TS_\ell v\leq s_\ell$. By the semidefinite constraint of \rf{f_G2} we have
\begin{align*}
v^TB\Delta&\leq v^T\left[\sum_\ell{\lambda}_\ell S_\ell\right]v+
\Delta^T\left[
\sum_k \mu_kT_k\right]
\Delta+
[A\Delta]^T\sum_i\gamma_ig_ig_i^TA\Delta\\
&\leq \sum_\ell{\lambda}_\ell s_\ell+
4\sum_k \mu_kt_k+\sum_i\gamma_i(g_i^TA\Delta)^2
\leq \phi_{\cS}(\lambda)+4\phi_\cT(\mu)+\sum_{i}\gamma_i(g_i^TA\Delta)^2\\
&\leq\phi_{\cS}(\lambda)+4\phi_\cT(\mu)+4\psi^2[G]\sum_{i}\gamma_i.
\end{align*}
Maximizing the left-hand side of the resulting inequality over $v\in\cB_*$, we arrive at $\|B\Delta\|\leq f_G(\lambda,\mu,\gamma)$. \qed
Note that the optimization problem in the right-hand side of (\ref{f_G2}) is an explicit convex optimization problem, so that $\Opt[G]$ is efficiently computable, provided that $\phi_\cS,\phi_\cT$ and $\pi$ are so. Thus, Proposition \ref{pro:G2} provides us with an efficiently computable upper bound on the $\epsilon$-risk of the polyhedral estimate stemming from a given contrast matrix $G$ and as such gives us a  computation-friendly tool to {\em analyse} the  performance of a polyhedral estimate. Unfortunately, this tool does not allow to {\em design} a ``presumably good'' estimate, since an attempt to make $G$ a variable, rather than a parameter, in the right-hand side problem in (\ref{f_G2}) results in a nonconvex, and thus difficult to solve, optimization problem.
We now look at two situations in which this difficulty can be overcome.

\subsection{Synthesis of ``presumably good'' contrast matrices}
We consider here two types of assumptions about the set $\cN$ of nuisances which allow for a computationally efficient design of ``presumably good'' contrast matrices. Namely,
\begin{enumerate}
\item[1)] ``ellitopic case:'' $\cN$ is a basic ellitope\\
\item[2)] ``co-ellitopic case:'' the set $N\cN=\{N\nu:\nu\in\cN\}$ is the {\sl polar} of the ellitope
\begin{align*}
&\cN_*=\{w\in\bR^m:\exists \overline{r}\in\overline{\cR}: w^T\ov R_jw\leq \overline{r}_j,j\leq \overline{J}\}\\
&\left[\hbox{ $\ov R_j\succeq0$, $\sum_j\ov R_j\succ0$; $\ov \cR\subset\bR^{\ov J}_+$, $\inter\ov \cR\neq\emptyset$, is a monotone convex compact}\right]\nonumber
\end{align*}
Note that $\cN_*$ is exactly the unit ball of the norm $\pi(g)=\max_{\nu\in \cN}g^TN\nu$.

\end{enumerate}
\subsubsection{Ellitopic case}
An immediate observation is that {\em the ellitopic case can be immediately reduced to the no-nuisance case}. Indeed, when $\cN$ is an ellitope, so is the direct product $\overline{\cX}=\cX\times\cN$. Thus, setting $\overline{A}[x;\nu]=Ax+N\nu$, $\ov{B}[x;\nu]=Bx$, observation (\ref{omega}) becomes
\[
\omega=\ov{A}\ov{x}_*+\xi,\eqno{[\overline{x_*}=[x_*;\nu_*]\in\overline{\cX}]}
\]
and our objective is to recover from this observation the linear image $w_*=\overline{B}\overline{x}_*$ of the new signal $\overline{x}_*$. Design of presumably good (and near-minimax-optimal when $\xi\sim\cN(0,\sigma^2I_m)$) polyhedral estimates in this setting is considered in \cite{PUP}. It makes sense to sketch the construction here, since it explains the idea which will be used through the rest of the paper.

Thus, consider the case when $\cN=\{0\}$, and let the signal set and the norm $\|\cdot\|$ still be given by (\ref{cXcB}). In this situation problem (\ref{f_G2}) becomes
\begin{align}\label{f_G2e}
\Opt[G]&=\min\limits_{\lambda,\mu,\gamma}\bigg\{
\phi_\cS(\lambda)+4\phi_\cT(\mu)+4\varkappa^2(\epsilon)\n^2[G]\sum_i\gamma_i:\;
\lambda\geq0,\;\mu\geq0,\,\gamma\geq0,\\
&\quad\quad\quad\quad
\left.\left[\begin{array}{c|c}\sum_\ell{\lambda}_\ell S_\ell&{1\over 2}B\\\hline
{1\over 2}B^T&
\sum_k \mu_kT_k+
A^T\left[\sum_i\gamma_ig_ig_i^T\right]A
\cr\end{array}\right]\succeq0\right\}\nonumber
\end{align}
Note that when $\theta>0$, we have $\Opt[G]=\Opt[\theta G]$. Indeed,  $(\lambda,\mu,\gamma)$ is a feasible solution to the problem specifying $\Opt[G]$
if and only if $\lambda,\mu,\theta^2\gamma)$ is a feasible solution to the problem specifying $\Opt[\theta G]$, and the values of the respective objectives
at these solutions are the same. It follows that as far as optimization of $\Opt[G]$ in $G$ is concerned, we lose nothing when restricting ourselves
to contrast matrices $G$ with $\varkappa(\epsilon)\n[G]=1$. In other words, setting
\be
\theta(g)=\varkappa(\epsilon)\|g\|_2
\ee{thetag} and  augmenting variables
$\lambda,\mu,$ and $\gamma$ in (\ref{f_G2e}) by variables $g_i$, $\theta(g_i)\leq1$,  $i=1,...,I$  (recall that we want to make $G$ variable rather than parameter and to minimize $\Opt[G]$ over $G$), we arrive at the problem
\begin{align}\label{f_G2f}
\Opt&=\min\limits_{\lambda,\mu,\gamma,\{g_i\},\rho}\bigg\{
\phi_\cS(\lambda)+4\phi_\cT(\mu)+4\rho:\;
\lambda\geq0,\;\mu\geq0,\,\gamma\geq0,\\
&\quad\quad
\left.\theta(g_i)\leq 1,\;\sum_i\gamma_i\leq \rho,\;\left[\begin{array}{c|c}\sum_\ell{\lambda}_\ell S_\ell&{1\over 2}B\\
\hline
{1\over 2}B^T&
\sum_k \mu_kT_k+
A^T\left[\sum_i\gamma_ig_ig_i^T\right]A
\cr\end{array}\right]\succeq0\right\}\nonumber
\end{align}
Now, aggregating variables $\gamma$, $g_1,...,g_I$ into the matrix $\Theta=\sum_i\gamma_ig_ig_i^T$ and denoting by $\mK$ the set of the pairs $(\Theta\in\bS^m_+,\rho)$ for which there exists decomposition $\Theta=\sum_{i\leq I}\gamma_ig_ig_i^T$ with $\theta(g_i)\leq 1$ and $\gamma_i\geq0$, $\sum_i\gamma_i\leq\rho$, (\ref{f_G2f}) can be rewritten as the optimization problem
\begin{align}\label{f_G2h}
\Opt&=\min\limits_{\lambda,\mu,\Theta,\rho}\bigg\{
\phi_\cS(\lambda)+4\phi_\cT(\mu)+4\rho:\;
\lambda\geq0,\;\mu\geq0,\\
&\quad\quad\left.(\Theta,\rho)\in\mK,\;
\left[\begin{array}{c|c}\sum_\ell{\lambda}_\ell S_\ell&{1\over 2}B\\
\hline
{1\over 2}B^T&
\sum_k \mu_kT_k+
A^T\Theta A
\cr\end{array}\right]\succeq0\right\}\nonumber
\end{align}
Observe that when $I\geq m$, $\mK$  is a simple convex cone:
\[\mK=\{(\Theta,\rho): \,\Theta\succeq0,\,\rho\geq \varkappa^2(\epsilon)\Tr(\Theta)\},
\]so that (\ref{f_G2h}) is an explicit (and clearly solvable) convex optimization program. To convert an optimal solution $(\lambda^*,\mu^*,\Theta^*,\rho^*)$ to \rf{f_G2h} into an optimal solution to \rf{f_G2f}, it suffices to subject $\Theta^*$ to the eigenvalue decomposition
$\Theta^*=\sum_{i=1}^I\upsilon_ie_ie_i^T$ with $\|e_i\|_2=1$ and $\upsilon_i\geq0$, $i\in\{1,..., m\}$, and $e_i=0$, $\upsilon_i=0$, $i\in\{m,...,I\}$, and set $g_i^*=\varkappa^{-1}(\epsilon)e_i$, $\gamma_i^*=\varkappa^2(\epsilon)\upsilon_i$, thus arriving  at an optimal solution $(\lambda^*,\,\mu^*,\,\{g_i^*,\gamma_i^*\}_{i\leq I},\,\rho^*)$ to problem (\ref{f_G2f}).

\subsubsection{Co-ellitopic case}
The just outlined approach to reducing the nonconvex problem (\ref{f_G2f}) responsible for the design of the best, in terms of $\Opt[G]$, contrast matrix $G$ to an explicit convex optimization problem heavily utilizes the fact that the unit ball of the norm $\theta(\cdot)$ (cf. \rf{thetag}) is a simple ellitope---a multiple of the unit Euclidean ball; this was the reason for $\mK$ to be a computationally tractable convex cone. Our future developments are built on
the fact that when the unit ball of
$\theta(\cdot)$ is a basic ellitope, something similar takes place: the associated set $\mK$, while not necessarily convex and computationally tractable, can be tightly approximated by a computationally tractable convex cone.  The underlying result (which is \cite[Proposition 3.2]{juditsky2024design}, up to notation) is as follows:
\begin{proposition}\label{propneed2} Let $I\geq m$, and let $\cW\subset\bR^m$ be a basic ellitope:
\begin{align*}
&\cW=\{w\in\bR^m:\,\exists r\in\cR:\,w^TR_jw\leq r_j,\,j\leq J\}\\
&\left[\hbox{$R_j\succeq0$, $\sum_jR_j\succ0;\;\cR\subset\bR^J_+$, $\inter\cR\neq\emptyset$, {\rm is a monotone convex compact}}\right]
\end{align*}
Let us associate with $\cW$ the closed convex cone\footnote{This indeed is a closed convex cone---the conic hull of the convex compact set $\{\Theta\succeq0:\; \exists  r\in\cR:\,\Tr(\Theta R_j)\leq r_j,\,1\leq j\leq J\}\times\{1\}$.}
$$
\bK=\left\{(\Theta,\rho):\, \exists r\in\cR:\,\Tr(\Theta R_j)\leq \rho r_j,\,1\leq j\leq J,\,\Theta\succeq 0,\,\rho\geq 0\right\}.
$$
Whenever a matrix $\Theta\in\bS^m_+$ is representable as $\sum_i\gamma_i w_iw_i^T$ with $\gamma_i\geq0$ and $w_i\in\cW$, one has $(\Theta,\sum_{i=1}^I\gamma_i)\in\bK$, and nearly vice versa:
whenever $(\Theta,\rho)\in\bK$,
one can find efficiently (via a randomized algorithm) vectors $w_i\in\cW$, and reals $\gamma_i\geq0$, $i\leq I$, such that $\Theta=\sum_i\gamma_iw_iw_i^T$ and
and
$$
\sum_i\gamma_i\leq 2\sqrt{2}\ln(4m^2J) \rho.
$$
\end{proposition}
We are now ready to outline a ``presumably good'' contrast design in the co-ellitopic case. Let us put $R_j=\four\ov R_j,\,j\leq \ov J$, and $R_{\ov J+1}=\tfrac{\varkappa^2(\epsilon)}{4}I_m$ and consider the ellitope
\begin{align}
\cW&=2\Big[\cN_*\cap \{w:\,\varkappa(\epsilon)\|w\|_2\leq1\}\Big]\nonumber\\
&=\{w\in\bR^m:\,\exists r\in\cR=\ov \cR\times[0,1]:\,w^TR_jw\leq r_j,\,j\leq J=\ov J+1\},\label{cW}
\end{align}
and let $\theta(\cdot)$ be the norm on $\bR^n$ with the unit ball $\cW$. Note that  $\theta(\cdot)=2\max\big[\pi(\cdot),\,\varkappa(\epsilon)\|\cdot\|_2\big]$, so that for every $G=[g_1,..,g_I]$, the quantity $\psi[G]$, see (\ref{f_G2}), is upper-bounded by $\max_i\theta(g_i)$, and this bound is tight within the factor 2. Consequently, Proposition \ref{pro:G2} states that the $\epsilon$-risk of the polyhedral estimate with contrast matrix $G$ is upper-bounded by the quantity
\begin{align}\label{f_G2i}
\ov\Opt[G]&=\min\limits_{\lambda,\mu,\gamma}\bigg\{
\phi_\cS(\lambda)+4\phi_\cT(\mu)+4[\max_i\theta(g_i)]^2\sum_i\gamma_i:\;
\lambda\geq0,\;\mu\geq0,\,\gamma\geq0,\\
&\quad\quad\quad\quad
\left.\left[\begin{array}{c|c}\sum_\ell{\lambda}_\ell S_\ell&{1\over 2}B\\\hline
{1\over 2}B^T&
\sum_k \mu_kT_k+
A^T\left[\sum_i\gamma_ig_ig_i^T\right]A
\cr\end{array}\right]\succeq0\right\}\nonumber
\end{align}
and $\Opt[G]\leq \ov\Opt[G]\leq 2\Opt[G]$. As in the previous section, the problem of minimizing $\ov\Opt[G]$ over $G$ can be reformulated in the form
(\ref{f_G2h}). A computationally efficient way to get a tight approximation to the optimal solution of the latter problem is given by the following result.
\begin{theorem}\label{thm:coell1}
Let $I\geq m$, and let
$$
\bK=\left\{(\Theta,\rho):\; \exists r\in\cR:\;\Tr(\Theta R_j)\leq \rho r_j,\,1\leq j\leq J,\,\Theta\succeq 0,\,\rho\geq 0\right\}
$$
(see (\ref{cW})).
Consider the convex optimization problem
\begin{align}\label{f_The2}
\Opt_*&=\min\limits_{\lambda,\mu,\gamma,\Theta,\rho}\bigg\{
\phi_\cS(\lambda)+4\phi_\cT(\mu)+4\underbrace{[2 \sqrt{2}\ln(4m^2J)]}_{=:\alpha}\rho
:\;
\lambda\geq0,\;\mu\geq0,\\
&\quad\quad\quad\quad\quad
\left.(\Theta,\rho)\in \bK,\;\left[\begin{array}{c|c}\sum_\ell{\lambda}_\ell S_\ell&{1\over 2}B\cr\hline
{1\over 2}B^T&
\sum_k \mu_kT_k+A^T
\Theta A
\cr\end{array}\right]\succeq0\right\}\nonumber
\end{align}
One can convert, in a computationally efficient way, the $\Theta$-component $\Theta^*$ of an optimal solution to this (clearly solvable) problem
into the contrast matrix $G^*$ such that
\[\ov\Opt[G^*]\leq {\sqrt{\alpha}}\min_G\ov\Opt[G]\leq 2 {\sqrt{\alpha}}\min_G\Opt[G].
\]
In particular, the $\epsilon$-risk of the polyhedral estimate with contrast matrix $G^*$ (this risk is upper-bounded by $\ov\Opt[G^*]$) does not exceed  $2\sqrt{\alpha}\min_G\Opt[G]$.
\end{theorem}
{\bf Proof.} When repeating the reasoning in the previous section, we conclude that $\overline{\Opt}:=\inf_G\ov\Opt[G]$ is equal to
\[
\inf_{g_1,...,g_I}\left\{\ov\Opt([g_1,...,g_I]):\;\max_i\theta(g_i)=1\right\}.
\] The latter $\inf$ is clearly attained at certain collection $g_i^+,...,g_I^+$ with $\max_i\theta(g_i^+)=1$. Let $G^+=[g_1^+,...,g_I^+]$, let $\lambda^+,\mu^+,\gamma_i^+,\,i\leq I,$ be an optimal solution to the problem in the right-hand side of (\ref{f_G2i}) associated with $g_i=g_i^+$, $i\leq I$, and let $\Theta^+=\sum_i\gamma_i^+[g_i^+][g_i^+]^T$ and $\rho^+=\sum_i\gamma_i^+$. We clearly have \[\ov\Opt=\ov\Opt[G^+]=\phi_\cS(\lambda^+)+4\phi_\cT(\mu^+)+4\rho^+.\]
Besides this, we are in the case where $\theta(g)\leq1$ is \aic{the same as}{equivalent to} $g\in\cW$, and therefore, by the first claim in
 Proposition \ref{propneed2}, $(\Theta^+,\rho^+)\in\bK$, implying that $(\lambda^+,\mu^+,\Theta^+,\rho^+)$ is a feasible solution to the optimization problem in (\ref{f_The2}).
 Due to the structure of the latter problem, for $\kappa>0$ the collection $(\kappa^{-1} \lambda^+,\kappa \mu^+,\kappa\Theta^+,\kappa\rho^+)$ is feasible for (\ref{f_The2}) with the corresponding value of the objective $\kappa^{-1}\phi_\cS(\lambda^+)+\kappa[\phi_\cT(\mu^+)+4\alpha\rho^+]$. It follows that
 \begin{align*}
 \Opt_*&\leq\inf_{\kappa>0}\left[\kappa^{-1}\phi_\cS(\lambda^+) +\kappa [4\phi_\cT(\mu^+) +4\alpha\rho^+]\right]\\
 &=2\big(\phi_\cS(\lambda^+)\underbrace{[4\phi_\cT(\mu^+)+4\alpha \rho^+]}_{\leq\alpha[4\phi_\cT(\mu^+)+4\rho^+]}\big)^{1/2}
 \leq2\sqrt{\phi_\cS(\lambda^+)[4\phi_\cT(\mu^+)+4\rho^+]} \sqrt{\alpha}\\&\leq \sqrt{\alpha}[\phi_\cS(\lambda^+)+4\phi_\cT(\mu^+)+4\rho^+]=\sqrt{\alpha}\,\overline{\Opt}.
 \end{align*}
 Finally, let $\ov\lambda,\ov\mu,\ov \Theta,\ov \rho$ be an optimal solution to (\ref{f_The2}).
 As $(\ov\Theta,\ov\rho)\in\bK$, the second claim in Proposition \ref{propneed2} states that there exists (and can be efficiently found) decomposition $\ov\Theta=\sum_i\ov\gamma_i[\ov g_i][\ov g]_i^T$ with $\ov g_i\in\cW$ (i.e., $\theta(\ov g_i)\leq1$), $i\leq I$, $\ov\gamma_i\geq0$, and $\sum_i\ov\gamma_i\leq\alpha \ov \rho$. The $\epsilon$-risk of the polyhedral estimate with the contrast matrix $\ov G=[\ov g_1,...,\ov g_I]$ is then upper-bounded by $\ov\Opt[\ov G]$. However, $\ov\lambda,\ov\mu,$ and $\{\ov\gamma_i\}$ form a feasible solution to the problem specifying $\ov\Opt[\ov G]$, and the value of the objective at this solution is upper bounded with
 \[\phi_\cS(\ov\lambda)+4\phi_\cT(\ov\mu)+4[\max_i\theta(\ov g_i)]\sum_i\ov\gamma_i
 \leq \phi_\cS(\ov\lambda)+4\phi_\cT(\ov\mu)+4\alpha\ov\rho=\Opt_*.
  \]Thus, the $\epsilon$-risk of the polyhedral estimate with contrast matrix $\ov G$ does not exceed
  \[\Opt_*\leq \sqrt{\alpha}\,\ov\Opt\leq 2\sqrt{\alpha}\min_G\Opt[G]. \eqno{\mbox{\qed}}\]

\section{Observations with outliers}\label{13estimate}
In this  section, we consider the estimation problem posed in Section \ref{theproblem} in the situation where the nuisance $\nu_*$ in (\ref{omega}) is sparse---has at most a given number $s$ of nonzero entries.
\paragraph{Estimate construction.}\label{theconstr}
Let  $\epsilon\in(0,1)$ be a given reliability tolerance. We consider the polyhedral estimate specified by two contrast matrices $H=[h_1,...,h_n]\in\bR^{m\times n}$ and $G=[g_1,...,g_I]\in\bR^{n\times I}$ which is as follows. Given observation $\omega$ (see (\ref{omega})) we solve the optimization problem
\beq\label{pro11}
\min_{\nu,x}\left\{\|\nu\|_1: \;x\in\cX,\;\begin{array}{ll}|h_k^T\left[N\nu+Ax-\omega\right]|\leq
\ov\varkappa(\epsilon)\|h_k\|_2,&k=1,...,n,\\
|g^T_i[N\nu+Ax-\omega]|\leq\ov\varkappa(\epsilon)\|g_i\|_2,&i=1,...,I,
\end{array}\right\}
\eeq
where
\[\ov\varkappa(\epsilon)=\sigma\sqrt{2\ln[2(n+I)/\epsilon]}.
\]
Let $(\widehat{\nu},\widehat{x})=(\widehat{\nu}(\omega),\widehat{x}(\omega))$ be
an optimal solution to the problem when the problem is feasible, otherwise we put $(\widehat{\nu},\widehat{x})=(0,0)$. Vector
\[
\wh w_{G,H}(\omega)=B\widehat{x}(\omega)
\] is the estimate of $w_*=Bx_*$ we intend to use.
\subsection{Risk analysis}\label{13analysis}
Let us denote $\Xi_\epsilon(G,H)$ the set of realizations of $\xi$ such that
\be|h^T_k\xi|\leq \ov\varkappa(\epsilon)\|h_k\|_2,\;k=1,...,n, \quad|g_i^T\xi|\leq \ov\varkappa(\epsilon)\|g_i\|_2,\;i=1,...,I,\quad\forall \xi\in \Xi_\epsilon(G,H).
\ee{chi11}
For the same reasons as in (\ref{neweq2}), one has
\[\Prob_{\xi\sim\SG(0,\sigma^2I_m)}(\Xi_\epsilon(G,H))\geq 1-\epsilon.
\]
Let us now fix $x_*\in \cX$, $s$-sparse $\nu_*$,  and $\xi\in\Xi_\epsilon(G,H)$, so that our observation is $\omega=Ax_*+N\nu_*+\xi$.
\paragraph{A.} By \rf{chi11} we have $|h_k^T\xi|\leq\ov\varkappa(\epsilon)\|h_k\|_2$ and $|g^T_i\xi|\leq\ov\varkappa(\epsilon)\|g_i\|_2$ for all $k\leq n$ and $i\leq I$,
while (\ref{pro11}) becomes the problem
\beq\label{becomes1}
\min_{{\nu},x}\left\{\|{\nu}\|_1:\;x\in\cX,\;\begin{array}{ll}|h^T_k\left[N[\nu-\nu_*]+A[x-x_*]-\xi\right]|\leq \ov\varkappa(\epsilon)\|h_k\|_2,&k=1,...,n,\\
|g^T_i\left[N[\nu-\nu_*]+A[x-x_*]-\xi\right]|\leq \ov\varkappa(\epsilon)\|g_i\|_2,&i=1,...,I.
\end{array}\right\}
\eeq
We conclude that $(\nu,x)=(\nu_*;x_*)$ is a feasible solution to (\ref{pro11}). Thus, we are in the case where $\widehat{\nu}$, $\widehat{x}$ are feasible for (\ref{becomes1}), and
\[
\|\widehat{\nu}\|_1\leq\|\nu_*\|_1.
\]
\paragraph{B.} Assume from now on that $(H,\|\cdot\|_\infty)$ satisfies Condition $Q_\infty(s,\kappa)$ of Section \ref{sectspars} with $\kappa<\half$,
that is,
\begin{equation}\label{ch1eq31p}
\|w\|_{\infty}\leq \|H^TNw\|_\infty+\tfrac{\kappa}s \|w\|_1\quad\forall w\in\bR^n.
\end{equation} Since $\widehat{\nu}$ and $\widehat{x}$ are feasible for \rf{becomes1}, we have
\[
|h_k^T[N[\widehat{\nu}-\nu_*]+A[\widehat{x}-x_*]-\xi]|\leq\ov\varkappa(\epsilon)\|h_k\|_2,\;\forall k\leq n.
\]
Invoking (\ref{chi11}) and the fact that $A[\wh x-x_*]\in2A\cX$ (since $\cX$ is symmetric w.r.t. the origin), we conclude that
\[
\|H^TN[\widehat{\nu}-\nu_*]\|_\infty\leq \max_{k}\left[\ov\varkappa(\epsilon)\|h_k\|_2+2\max_{x\in\cX}|h_k^TAx|\right],
\]
and besides this, $\nu_*$ is $s$-sparse and $\|\widehat{\nu}\|_1\leq\|\nu_*\|_1$.
Now  Proposition \ref{newl1theorem} with $\nu_*$ in the role of ${\nu}$ 
implies that
\be
\|\widehat{\nu}-\nu_*\|_q
\leq \tfrac{(2s)^{\frac1q}}{1-2\kappa}\max_{k}\left[\ov\varkappa(\epsilon)\|h_k\|_2+2\max_{x\in\cX}|h_k^TAx\|\right],\quad1\leq q\leq \infty,
\ee{ql1}
in particular, that
\begin{subequations}
\label{quick10}
 \begin{align}\label{quick10a}
  \|\widehat{\nu}-\nu_*\|_\infty&\leq{1\over 1-2\kappa}\max_{k}\left[\ov\varkappa(\epsilon)\|h_k\|_2+2\max_{x\in\cX}|h_k^TAx\|\right]
  =:\rho_H,\\
\|\widehat{\nu}-\nu_*\|_1&\leq{2s\over 1-2\kappa}\max_{k}\left[\ov\varkappa(\epsilon)\|h_k\|_2+2\max_{x\in\cX}|h_k^TAx\|\right]
=2s\rho_H.
 \label{quick10b}
 \end{align}
\end{subequations}
In addition,  \cite[Proposition 1.3.4]{PUP} states that the set $\cH$ of the pairs $(H,\kappa)$ with $m\times n$ matrices $H$ satisfying {Condition $Q_\infty(s,\kappa)$} is the computationally tractable convex set
\beq\label{quick4}
\cH=\left\{(H,\kappa)\in\bR^{m\times n}\times\bR:\;\big|[I_n-N^TH]_{ij}\big|\leq s^{-1}\kappa,\,1\leq i,j\leq n\right\}.
\eeq
\paragraph{C.} Since $\widehat{\nu}$ and $\widehat{x}$ are feasible for (\ref{becomes1}), we have
$$
|g_i^T\left(N[\widehat{\nu}-\nu_*]+A[\widehat{x}-x_*]-\xi\right)|\leq\ov\varkappa(\epsilon)\|g_i\|_2,\quad i=1,...,I,
$$
while $|g_i^T\xi|\leq  \ov\varkappa(\epsilon)\|g_i\|_2 \;\forall i$ due to
$\xi\in\Xi_\epsilon(G,H)$. We conclude that
\begin{equation}\label{deq12}
|g_i^TA[\widehat{x}-x_*]|\leq \aic{}{2}\ov\varkappa(\epsilon)\|g_i\|_2+|g^T_iN[\widehat{\nu}-\nu_*]|,\;i\leq I.
\end{equation}
Let $\|z\|_{k,1}$, $z\in\bR^n$, be the sum of $\min[k,n]$ largest magnitudes of entries in $z$; note that $\|\cdot\|_{k,1}$ is the norm conjugate to the norm with the unit ball $\{u:\|u\|_\infty\leq1,\|u\|_1\leq k\}$. Consequently,
 (\ref{quick10}) implies that
\begin{equation}\label{concl}
|g_i^TN[\widehat{\nu}-\nu_*]|\leq \rho_H \|N^Tg_i\|_{2s,1},
\end{equation} and, therefore, by (\ref{deq12})
\beq\label{quick11}
|g_i^TA[\widehat{x}-x_*]|\leq \psi_H[G], \,\,\psi_H[G]=\max_i\big[\aic{}{2}
\ov\varkappa(\epsilon)\|g_i\|_2+\rho_H \|N^Tg_i\|_{2s,1}\big].
\eeq
Applying the same argument as in the proof of Proposition \ref{pro:G2}, with (\ref{quick11}) in the role of (\ref{role}), we arrive at the following result:
\begin{proposition}\label{propro} In the situation of this section given $\kappa\in(0,1/2)$ and $m\times n$ matrix $H$ satisfying $(H,\kappa)\in\cH$, see (\ref{quick4}),
consider optimization problem (cf. \rf{f_G2})
\begin{align}\label{f_G2+}
\Opt[G,H]&=\min\limits_{\lambda,\mu,\gamma}\bigg\{
\overbrace{\phi_\cS(\lambda)+4\phi_\cT(\mu)+
\aic{\left[\max_i[\ov\varkappa(\epsilon)\|g_i\|_2+\rho_H\|N^Tg_i\|_{2s,1}]\right]^2}{\psi^2_H[G]}\sum_i\gamma_i}^{=:f_{G,H}(\lambda,\mu,\gamma)}:\nn
&\lambda\geq0,\;\mu\geq0,\,\gamma\geq0,
\left.\left[\begin{array}{c|c}\sum_\ell{\lambda}_\ell S_\ell&{1\over 2}B\\\hline
{1\over 2}B^T&
\sum_k \mu_kT_k+
A^T\left[\sum_i\gamma_ig_ig_i^T\right]A
\cr\end{array}\right]\succeq0.\right\}
\end{align}
Let $(\lambda,\mu,\gamma)$ be a feasible solution to \rf{f_G2+}. Then
\[\Risk_{\epsilon}[\widehat{w}_{G,H}\aic{|\cX,\cN}{}]\leq f_{G,H}(\lambda,\mu,\gamma).
\]\end{proposition}
\subsection{Synthesis of contrast matrices}\label{pgc1}
Our present objective is to design contrast matrices $H$ and $G$ with small value of the bound $\Opt[G,H]$ for the $\epsilon$-risk of the estimate $\wh w_{G,H}$.
\paragraph{D.} Building the contrast matrix $H\in \bR^{m\times n}$ is straightforward: the risk bound $\Opt[G,H]$, depends on $H=[h_1,...,h_n]$ solely through the quantity
\def\err{\hbox{\scriptsize RecoveryError}}
\[
\rho_H={1\over 1-2\kappa}\max_{k\leq n}\left[\ov\varkappa(\epsilon)\|h_k\|_2+2\max_{x\in \cX}\|h_k^TAx\|_\infty\right].
\]
and is smaller the smaller is $\rho_H$. For a fixed $\kappa\in(0,1/2)$, a presumably good choice of $H=[h_1,...,h_n]$ is then given by  optimal solutions to $n$ optimization problems
\be
h_k=\argmin_{h}\bigg\{\ov\varkappa(\epsilon)\|h\|_2+2\max_{x\in\cX}\|h^TAx\|_\infty:\; h\in\bR^{m},\,\|\Col_i[I_n-N^Th]\|_\infty\leq s^{-1}\kappa\bigg\}
\ee{hko}
which, when recalling what $\cX$ is, by conic duality, are equivalent to the problems
\begin{align*}
h_k&=\argmin_{h,v,\chi}\bigg\{\ov\varkappa(\epsilon)\|h\|_2+ v+\phi_{\cT}(\chi):\; h\in\bR^{m},\chi\geq 0,\nn
&\quad \quad\quad \quad\quad \left[\begin{array}{c|c}v&h^TA\\\hline
A^Th&\sum_k\chi_k T_k\end{array}\right]\succeq0,\;\|\Col_i[I_n-N^Th]\|_\infty\leq s^{-1}\kappa\bigg\},\quad 1\leq k\leq n.
\end{align*}
\paragraph{E.}
The proposed construction of $G$ is less straightforward.
We  proceed as follows. Let $G=[G_1,G_2]$ where $G_2, G_1\in \bR^{m\times m}$ (so that $I=2m$).
\paragraph{E.1}
Notice that as $\xi\in\Xi_\epsilon(G,H)$, problem (\ref{becomes1}) is feasible, and  $(\wh x,\wh\nu)$ is its feasible solution. For a column $g$ of $G$, by the constraints of the problem, we have
\begin{equation}\label{deq13}
|g^TA[\widehat{x}-x_*]|\leq \aic{}{2}\ov\varkappa(\epsilon)\|g\|_2+|g^TN[\widehat{\nu}-\nu_*]|
\leq \aic{}{2}\ov\varkappa(\epsilon)\|g\|_2+2s\rho_H\|N^Tg\|_\infty,
\end{equation}
(we have used  (\ref{deq12}) and (\ref{concl})),
implying that
\be
\big(g^TA[\widehat{x}-x_*]\big)^2\leq 2\left(\aic{}{4}\ov\varkappa^2(\epsilon)\|g\|_2^2+4s^2\rho^2_H\|N^Tg\|^2_\infty\right),\quad i=1,...,m.
\ee{g1cone}
Note  that the set
\[
\cM=\left\{g\in \bR^m: \,\aic{2}{8}\ov\varkappa^2(\epsilon)\|g\|_2^2+8s^2\rho^2_H\|N^Tg\|^2_\infty\leq 1\right\}
\]
is an ellitope: when denoting $N=[\n_1,...,\n_n]$ we have
\[
\cM=\bigg\{g\in \bR^m: \;\exists r\in[0,1]^n:\;g^T\underbrace{\left(\aic{2}{8}\ov\varkappa^2(\epsilon)I_m+8s^2\rho^2_H\n_j\n_j^T\right)}_{M_j}g\leq r_j,\,j=1,...,n\bigg\}.
\]
\paragraph{E.2} Next, observe that when $\xi\in \Xi_\epsilon(G,H)$, by \rf{ql1} one has
\[\|\widehat{\nu}-\nu_*\|_2
\leq \tfrac{\sqrt{2s}}{1-2\kappa}\max_{k\leq n}\left[\ov\varkappa(\epsilon)\|h_k\|_2+\max_{x\in\cX}|h_k^TAx|\right]=\sqrt{2s}\rho_H.
\]
Then by (\ref{deq13}), for a column $g$ of $G$ it holds
\begin{align}
\left(g^TA[\widehat{x}-x_*]\right)^2&\leq \left(\aic{}{2}\ov\varkappa(\epsilon)\|g\|_2+|g^TN[\widehat{\nu}-\nu_*]|\right)^2\leq \left(
\aic{}{2}\ov\varkappa(\epsilon)\|g\|_2+\sqrt{2s}\rho_H \|N^Tg\|_{2}\right)^2\nn
&\leq g^T\left(\aic{2}{8}
\ov\varkappa^2(\epsilon)I_m+4s\rho^2_H NN^T\right)g.\label{qt201}
\end{align}
Now, let us put
\be
Q=\big(\aic{2}{8}
\ov\varkappa^2(\epsilon)I_m+4s\rho^2_H NN^T\big)^{-1/2},
\ee{qrf}
and consider the optimization problem
\begin{subequations}\label{33old}
\begin{align}\label{f_The22}
\Opt&=\min\limits_{\lambda,\mu,\Theta_1,\Theta_2,\rho}\bigg\{f_{H}(\lambda,\mu,\Theta_1,\Theta_2,\rho):\;
\lambda\geq0,\;\mu\geq0,\,\Theta_1\succeq0,\,\Theta_2\succeq0,\\
&\quad\quad\quad\Tr(M_j\Theta_1)\leq \rho,\;j=1,...,n,\;
\left.\left[\begin{array}{c|c}\sum_\ell{\lambda}_\ell S_\ell&{1\over 2}B\cr\hline
{1\over 2}B^T&
\sum_k \mu_kT_k+A^T
(\Theta_1+Q\Theta_2Q^T) A
\cr\end{array}\right]\succeq0\right\}\nonumber
\end{align}
where
\begin{align}\label{f_The2f}
f_{H}(\lambda,\mu,\Theta_1,\Theta_2,\rho)&=\phi_\cS(\lambda)+4\phi_\cT(\mu)+\Tr(\Theta_2)
+2\sqrt{2}\ln(4m^2n)\rho.
\end{align}
\end{subequations}
Note that the constraints on $\Theta_1$ and $\rho$ of the problem \rf{f_The22} say exactly that $(\Theta_1,\rho)$ belongs to the cone $\bK$ associated, as explained in Proposition \ref{propneed2},  with the ellitope $\cM$ in the role of $\cW$.
\begin{theorem}\label{the:newsparse} Given a feasible solution
$(\lambda,\mu, \tau,\Theta_1, \Theta_2)$ to \rf{33old}, let us build $m\times m
$ contrast matrices  $G_1$, $G_2$ as follows.
\begin{itemize}
\item To build $G_1$, we apply the second part of Proposition \ref{propneed2} to $\Theta_1,\rho,\cM$ in the roles of $\Theta,\rho,\cW$, to get, in a computationally efficient way, a decomposition $\Theta_1=\sum_{i=1}^m\gamma_ig_{1,i}g_{1,i}^T$ with $g_{1,i}\in\cM$ and $\gamma_i\geq0$, $\sum_i\gamma_i\leq 2\sqrt{2}\ln(4m^2n)\rho$. We set $G_1=[g_{1,1},...,g_{1,m}]$.
\item To build $G_2$, we subject $\Theta_2$ to eigenvalue decomposition $\Theta_2=\Gamma\Diag\{\chi\}\Gamma^T$ and set $G_2=[g_{2,1},...,g_{2,m}]=Q\Gamma$.
\end{itemize}
Note that $\Theta_1+Q\Theta_2Q=\sum_i\gamma_ig_{1,i}g_{1,i}^T+\sum_i\chi_ig_{2,i}g_{2,i}^T$.

For the resulting polyhedral estimate $\wh w_{G,H}$ and for all $x_*\in\cX$, $s$-sparse $\nu_*$ and $\xi\in\Xi_\epsilon(G,H)$ it holds
\begin{equation}\label{target}
\|\wh w_{G,H}(Ax_*+N\nu_*+\xi)-Bx_*\|\leq f_H(\lambda,\mu,\Theta_1,\Theta_2,\rho)
\end{equation}
implying that the $\epsilon$-risk of the estimate is upper-bounded by $f_H(\lambda,\mu,\Theta_1,\Theta_2,\rho)$ (due to $\xi\in\Xi_\epsilon(G,H)$
with probability $\geq1-\epsilon$).
\end{theorem}
{\bf Proof.} Let us fix $x_*\in\cX$, $s$-sparse $\nu_*$, $\xi\in\Xi_\epsilon(G,H)$, and let $w=Ax_*+N\nu_*+\xi$. By {\bf A},
problem  (\ref{pro11}) is feasible, so that $(\wh x,\wh \nu)=(\wh x(\omega),\wh \nu(\omega))$ is its optimal solution, and $\wh w=B\wh x$ is the estimate $\wh w_{G,H}(\omega)$. Let  $\Delta=\wh x-x_*$, and let $e_1,...,e_m$ be the columns of the orthonormal matrix $\Gamma$. By construction of $G_2$, we have for all $ j\leq m$ (see \rf{qt201})
 \begin{align*}
(g_{2,j}^TA\Delta)^2&\leq  g_{2,j}^T \left(\aic{2}{8}
\ov\varkappa^2(\epsilon)I_m+4s\rho^2_H NN^T\right)g_{2,j}=e_j^T[Q\left(\aic{2}{8}
\ov\varkappa^2(\epsilon)I_m+4s\rho^2_H NN^T\right)Q]e_j
=e_j^Te_j=1.
\end{align*}
 Furthermore, due to $g_{1,i}\in\cM$ one has (see \rf{g1cone})
 \[
 \big(g^T_{1,i}A\Delta\big)^2\leq \aic{2}{8}\ov\varkappa^2(\epsilon)\|g\|_2^2+8s^2\rho^2_H\|N^Tg\|^2_\infty\leq 1\quad \forall i\leq m.
 \]
Now, by the semidefinite constraint of \rf{f_The22} and due to $\Theta_1+Q\Theta_2Q=\sum_i\gamma_ig_{1,i}g_{1,i}^T+\sum_i\chi_ig_{2,i}g_{2,i}^T$, for every $v\in \cB_*$ we have
\begin{align*}
v^TB\Delta&\leq v^T\left[\sum_\ell{\lambda}_\ell S_\ell\right]v+
\Delta^T\left[
\sum_k \mu_kT_k\right]
\Delta+
[A\Delta]^T\left[\sum_i\gamma_ig_{1,i}g_{1,i}^T+\sum_i\chi_ig_{2,i}g_{2,i}^T\right]A\Delta\\
&\leq \phi_{\cS}(\lambda)+4\phi_\cT(\mu)+\sum_{i}\chi_i(g_{1,i}^TA\Delta)^2+\sum_{j}\gamma_j(g_{2,j}^TA\Delta)^2\\
&\left[\hbox{as $[v^TS_1v;...;v^TS_Lv]\in\cS$ due to $v\in\cB_*$ and $[\Delta^TT_1\Delta;...;\Delta^TT_L\Delta]\in 4\cT$  due to $\Delta\in2\cX$}\right]\\
&\leq \phi_{\cS}(\lambda)+4\phi_\cT(\mu)+\sum_{i}\chi_i+\sum_{j}\gamma_j\leq f_{H}(\lambda,\mu,\tau,\Theta_1,\Theta_2)
\end{align*}
due to $\sum_i\gamma_i\leq 2\sqrt{2}\ln(4m^2n)\rho$ and $\sum_i\chi_i=\Tr(\Theta_2)$.
Taking the supremum over $v\in\cB_*$ in the resulting inequality, we arrive at (\ref{target}). \qed

\subsection{An alternative}\label{sectaltern}
\aic{}{Our objective in this section is to refine risk bounds \rf{f_G2+} and \rf{f_The22} to produce more efficient contrasts. Our course of action is as follows. First, to extend the possible choice of $H$-contrasts`` responsible'' for the perturbation recovery, we refine the bounds \rf{quick10} for the accuracy of sparse recovery, notably,  to allow using contrasts not satisfying Condition $Q_{\infty}(s,\kappa)$. Second, we improve the bounding of the risk of the estimate $\wh w(\omega)$ by taking into account the contribution of the $H$-component of the ``complete'' contrast matrix $[H,G]$ when optimizing the $G$-component of the contrast.}

In the sequel, we consider the estimate described in the beginning of Section \ref{13estimate},
the only difference being in the sizes of contrast matrices $G$ and $H$: now $H=[h_1,...,h_M]\in \bR^{m\times M}$, and $G=[g_1,...,g_{2m}]$. Thus, in our present setting, given observation $\omega$, we solve the optimization problem
\be
\min_{{\nu},x}\left\{\|{\nu}\|_1:\,x\in\cX,\,\begin{array}{ll}|{h^T_k(N\nu}+Ax-\omega)|\leq\ov\varkappa(\epsilon)\|h_k\|_2,&k=1,...,M,\\
|g_{i}^T({N\nu}+Ax-\omega)|\leq \ov\varkappa(\epsilon)\|g_{2,i}\|_2,&i=1,...,2m,\end{array}\right\}
\ee{pro1a}
with
\[\ov\varkappa(\epsilon)=\sigma\sqrt{2\ln[(2M+4m)/\epsilon]}\],
specify $\wh x(\omega),\wh \nu(\omega)$ as an optimal solution to the problem when the problem is feasible, otherwise set $(\wh x(\omega),\wh \nu(\omega))=(0,0)$, and take $\wh w_{G,H}(\omega)=B\wh x(\omega)$ as the estimate of $Bx_*$.
\subsubsection{Risk analysis}\label{modest}
The above problem can be rewritten equivalently as
\be
\min_{{\nu},x}\left\{\|{\nu}\|_1:\,x\in\cX,\,\begin{array}{ll}|{h^T_k(N[\nu-\nu_*]}+A[x-x_*]-\xi)|\leq\ov\varkappa(\epsilon)\|h_k\|_2,&k=1,...,M,\\
|g_{i}^T({N[\nu-\nu_*]}+A[x-x_*]-\xi)|\leq \ov\varkappa(\epsilon)\|g_{i}\|_2,&i=1,...,2m,
\end{array}\right\}
\ee{problem14}
and when setting
\be
\Xi_\epsilon(G,H):=\left\{\xi\in\bR^m:\,\begin{array}{ll}
|h_k^T\xi|\leq \ov\varkappa(\epsilon)\|h_k\|_2,&k=1,...,M,\\
|g_{i}^T\xi|\leq \ov\varkappa(\epsilon)\|g_{i}\|_2,&i=1,...,2m,\end{array}\right\}
\ee{Xi}
we have
$$
\Prob_{\xi\sim\SG(0,\sigma^2I_m)}({\Xi_\epsilon(G,H)})\geq 1-\epsilon.
$$
Let us fix $\xi\in{\Xi_\epsilon(G,H)}$  and set $\omega=Ax_*+N\nu_*+\xi$. As
$(\widehat{\nu},\,\widehat{x})$ is a feasible for (\ref{problem14}), $\wh x:=\wh x(\omega)$, $\wh\nu:=\wh\nu(\omega)$ is feasible as well, $\|\widehat{\nu}\|_1\leq\|\nu_*\|_1$. Thus, same as in the proof of Proposition \ref{newl1theorem}, for $z=\wh \nu-\nu_*$ it holds
\[\|z\|_1\leq 2\|z\|_{s,1}\]
implying that
\begin{equation}\label{wehave}
\|z\|_1\leq 2s\|z\|_\infty,\quad \|z\|_2\leq \sqrt{2s}\|z\|_\infty.
\end{equation}
Now denote
$
\Delta=\widehat{x}-x_*,
$
and consider $n$ pairs of convex optimization problems
\begin{align}
\Opt_2[i]&=\max_{v,t,w}\left\{\sqrt{w_it}: v\in 2\cX,\,\begin{array}{l}
\|w\|_\infty\leq w_i,\;\|w\|_1\leq t,\;t\leq 2sw_i,\\
|h_k^T(Nw+Av)|\leq 2\ov\varkappa(\epsilon)\|h_k\|_2,\,k=1,...,M\end{array}\right\}\tag{$P_2[i]$}\label{p2}\\
\Opt_\infty[i]&=\max_{v,w}\left\{w_i: v\in 2\cX,\,\begin{array}l
\|w\|_\infty\leq w_i,\,\|w\|_1\leq 2sw_i,\\
|h_k^T(Nw+Av)|\leq 2\ov\varkappa(\epsilon)\|h_k\|_2,\,k=1,...,M.\end{array}\right\}&\tag{$P_\infty[i]$}\label{pinf}
\end{align}\noindent
Observe that a feasible solution $(v,t,w)$ to $(P_2[i])$ satisfies $\|w\|_\infty\leq w_i$ and $\|w\|_1\leq t$, whence
\be\|w\|_2\leq \sqrt{w_it}\leq\Opt_2[i].\ee{Opt2i!}
Now, let $\iota=\iota_z$ be the index of the largest in magnitude entry in $z$. Taking into account that $\xi\in{\Xi_\epsilon(G,H)}$
and recalling that $\Delta\in2{\cX}$, we conclude that when $z_\iota\geq0$, $(v,t,w)=(\Delta,\|z\|_1,z)$ is feasible for $(P_2[\iota])$
and $(v,w)=(\Delta,z)$ is feasible for $(P_\infty[\iota])$, while when $z_\iota<0$ the same holds true for $(v,t,w)=(-\Delta,\|z\|_1,-z)$ and
$(v,w)=(-\Delta,-z)$.
Indeed in the first case $v=\Delta\in \cX$, $|h_k^T[A\wh x+N\wh \nu-\omega]|\leq\ov\varkappa(\epsilon)\|h_k\|_2$ and
$|h_k^T[Ax_*+N\nu_*-\omega]|\leq\ov\varkappa(\epsilon)\|h_k\|_2$ as both pairs $(\wh x,\wh \nu)$ and $(x_*,\nu_*)$ are feasible for (\ref{pro1a}), implying the second line constraints of \rf{p2}.  Note that we are in the case of $z_\iota=\|z\|_\infty$, that is, constraints in the first line of \rf{p2} are satisfied for $w=z$ due to (\ref{wehave}).
Thus, $(\Delta,\|z\|_1,z)$ indeed is feasible for (\ref{p2}). As a byproduct of our reasoning, $(\Delta,z)$ is feasible
for (\ref{pinf}).  In the second case, the reasoning is completely similar.
\par
Next, setting
\be
\Opt_2=\max_i\Opt_2[i],\quad\Opt_\infty=\max_i\Opt_\infty[i],
\ee{Opt2max}
and recalling that $(\Delta,\|z\|_1,z)$ or $(-\Delta,\|z\|_1,-z)$ is feasible for some of the problems $(P_2[i])$, and $(\Delta,z)$ or $(-\Delta,-z)$ is feasible for some of the problems $(P_\infty[i])$, when invoking \rf{Opt2i!} we get
for all $\xi\in{\Xi_\epsilon(G,H)}$
\[
\|z\|_\infty\leq\Opt_\infty,\quad\|z\|_2\leq\Opt_2,\quad
\|z\|_1\leq 2s\Opt_\infty.
\]
Consequently, for all $d\in \bR^m$
\begin{align}
|d^TNz|&\leq \max_{z}\left\{d^TNz:\,\|z\|_\infty\leq\Opt_\infty,\,\|z\|_2\leq\Opt_2,\,
\|z\|_1\leq 2s\Opt_\infty\right\}\nn
&=\underbrace{\min_{u,v,w}\left\{\|u\|_1\Opt_\infty+\|v\|_2\Opt_2+2s\|w\|_\infty\Opt_\infty,\,u+v+w=N^Td\right\}}_{=:\ov\pi(N^Td)}.\label{eq41}
\end{align}
Now, recalling that $\wh x,\wh\nu$ is feasible for (\ref{problem14}) and that $\xi\in\Xi_\epsilon(G,H)$, we conclude that \aic{}{columns $d_i$, $i=1,...,M+2m$ of the {\em ``aggregated'' contrast matrix} $D=[G,H]$ satisfy}
\[|\aic{g_i}{d_i}^TA\Delta|\leq |\aic{g_i}{d_i}^TNz|+|\aic{g_i}{d_i}^T\xi|+\ov\varkappa(\epsilon)\|g\|_2 ,
\] whence
\begin{equation}\label{latterr}
|\aic{g_i}{d_i}^TA\Delta|\leq \underbrace{\ov\pi(N^T\aic{g_i}{d_i})+2\ov\varkappa(\epsilon)\|\aic{g_i}{d_i}\|_2}_{\aic{}{=:\psi_H(d_i)}},\quad i\leq \aic{}{M+}2m.
\end{equation}
Applying the same argument as in the proof of Proposition \ref{pro:G2}, with (\ref{latterr}) in the role of (\ref{role}), we arrive at the following result:
\begin{proposition}\label{proproplus} In the situation of this section,
consider optimization problem (cf. \rf{f_G2+})
\begin{align}\label{f_G2++}
\Opt[G,H]&=\min\limits_{\lambda,\mu,\gamma}\bigg\{\overbrace{\phi_\cS(\lambda)+4\phi_\cT(\mu)+\sum_i\gamma_i
\aic{\left[\ov\pi(N^T\aic{g_i}{d_i})
+2\ov\varkappa(\epsilon)\|\aic{g_i}{d_i}\|_2\right]^2}{\psi^2_H(d_i)}}^{=:\bar f_{G,H}(\lambda,\mu,\gamma)}:\,
\lambda\geq0,\;\mu\geq0,\,\gamma\geq0,\nn
&\quad\quad\quad\quad
\left.\left[\begin{array}{c|c}\sum_\ell{\lambda}_\ell S_\ell&{1\over 2}B\\\hline
{1\over 2}B^T&
\sum_k \mu_kT_k+
A^T\left[\sum_i\gamma_i\aic{g_i}{d_i}\aic{g_i}{d_i}^T\right]A
\cr\end{array}\right]\succeq0.\right\}
\end{align}
Let $(\lambda,\mu,\gamma)$ be a feasible solution to \rf{f_G2++}. Then
\[\Risk_{\epsilon}[\widehat{w}_{G,H}|\cX,\cN]\leq \bar f_{G,H}(\lambda,\mu,\gamma).
\]\end{proposition}
\subsubsection{Contrast matrix synthesis}\label{synt}
We continue our analysis of the estimate $\wh w_{G,H}$ in the situation when the observation is $\omega=Ax_*+N\nu_*+\xi$ with $\xi\in\Xi_\epsilon(G,H)$, see (\ref{Xi}). By (\ref{eq41}), for $z=\wh\nu-\nu_*$ and all $g\in\bR^m$ we have
\begin{align*}
|g^TNz|\leq\min\left\{\|N^Tg\|_2 \Opt_2,\;\sqrt{2s}\|N^Tg\|_2 \Opt_\infty,\;2s\|N^Tg\|_\infty \Opt_\infty
\right\}
\end{align*}
what implies (cf. \rf{latterr}) that for all $i\leq 2m$
\be
|g^T_iA\Delta|\leq 2\ov\varkappa(\epsilon)\|g_i\|_2+\min\left\{\|N^Tg_i\|_2 \Opt_2,\;\sqrt{2s}\|N^Tg_i\|_2 \Opt_\infty,\;2s\|N^Tg_i\|_\infty \Opt_\infty
\right\}.
\ee{gTAD1}
Note that the right-hand side in (\ref{gTAD1}) is nonconvex in $g$, making our design techniques unapplicable.
To circumvent this difficulty, we intend to utilize the following important feature of polyhedral estimates: one may easily ``aggregate'' several estimates of this type to yield an estimate with the risk which is nearly as small as the smallest of the risks of the estimates combined.

Here is how it works in the present setting. We split the $m\times 2m$ contrast $G$ into two $m\times m$ blocks $G_\chi=[g_{\chi,1},...,g_{\chi,m}]$, $\chi=1,2$, and design the blocks utilizing the respective inequalities inherited from (\ref{gTAD1}), specifically, the inequalities
\begin{align*}
|g^T_{1,i}A\Delta|&
\leq 2\ov\varkappa(\epsilon)\|g_{1,i}\|_2+2s\|N^Tg_{1,i}\|_\infty \Opt_\infty,\\
|g^T_{2,i}A\Delta|&
\leq2\ov\varkappa(\epsilon)\|g_{2,i}\|_2+\|N^Tg_{2,i}\|_2\underbrace{\min\{\Opt_2,\,\sqrt{2s}\Opt_\infty\}}_{=:\varrho_{2,H}}
\end{align*}
where $\Delta=\wh x-x_*$.
We weaken these inequalities to
\begin{align*}
|g^T_{1,i}A\Delta|^2&\leq\pi^2_1(g_{1,i}),\quad\pi_1(g)=\sqrt{8\ov\varkappa^2(\epsilon)\|g\|_2^2+8s^2\Opt_\infty^2\|N^Tg\|_\infty^2},\\
|g^T_{2,i}A\Delta|^2&\leq \pi_2^2(g_{2,i}),\quad\pi_2(g)=\sqrt{8\ov\varkappa^2(\epsilon)\|g\|_2^2+2\varrho_{2,H}^2\|N^Tg\|_2^2}.
\end{align*}
Notice that norms $\pi_\chi$, $\chi=1,2$, are ellitopic, so we can use in our present situation the techniques from Section \ref{pgc1}, thus arriving at an analogue of Theorem \ref{the:newsparse}. To this end, denote by $\n_1,...,\n_n$ the columns of $N$ and set
\[
\ov M_j=8\ov\varkappa^2(\epsilon)I_m+8s^2\Opt_\infty^2\n_j\n_j^T,\;j\leq m,\;\;\mbox{and}\;\;
\ov Q=\left(8\ov\varkappa^2(\epsilon)I_m+2\varrho_{2,H}^2NN^T\right)^{-1/2}.
\]
Next, observe that the unit ball of $\pi_1(\cdot)$ is the ellitope
\[
\ov \cM=\{w\in\bR^m:\exists r\in[0,1]^M: w^T\ov M_j w\leq \rho_j,\,j=1,...,M\}
\]
and the unit ball of $\pi_2$ is the ellipsoid $w^T\ov Q^{-2}w \leq 1$. Now, let us consider the optimization problem
\begin{subequations}\label{45old}
\begin{align}\label{f_The22o}
\Opt=&\min\limits_{\lambda,\mu,\tau,\Theta_1,\Theta_2,\rho}\bigg\{\ov f_{H}(\lambda,\mu,\aic{}{\tau,}\Theta_1,\Theta_2,\rho):\;
\lambda\geq0,\;\mu\geq0,\,\aic{}{\tau\geq 0,}
\\& \left.\begin{array}{l}
\Theta_1\succeq0,\,\Theta_2\succeq0,\;\Tr(\ov M_j\Theta_1)\leq \rho,\;j=1,...,n,\\
\left[\begin{array}{c|c}\sum_\ell{\lambda}_\ell S_\ell&{1\over 2}B\cr\hline
{1\over 2}B^T&
\sum_k \mu_kT_k\aic{}{+A^T\left[\sum_i\tau_ih_ih_i^T\right]A}+A^T
(\Theta_1+Q\Theta_2Q^T) A
\cr\end{array}\right]\succeq0\end{array}\right\}\nonumber
\end{align}
where
\begin{align}\label{f_The2fo}
\ov f_{H}(\lambda,\mu,\tau,\Theta_1,\Theta_2,\rho)&=\phi_\cS(\lambda)+4\phi_\cT(\mu)\aic{}{+\sum_i\tau_i\psi^2_H(h_i)}+\Tr(\Theta_2)
+2\sqrt{2}\ln(4m^2n)\rho.
\end{align}
\end{subequations}
Note that the constraints on $\Theta_1$ and $\rho$ in this problem say exactly that $(\Theta_1,\rho)$ belongs to the cone $\bK$ associated, according to Proposition \ref{propneed2},  with the ellitope $\ov\cM$ in the role of $\cW$.
\begin{theorem}\label{the:newsparseo} Given a feasible solution
$(\lambda,\mu, \aic{}{\tau,}\Theta_1, \Theta_2,\rho)$ to \rf{45old}, let us build $m\times m
$ contrast matrices  $G_1$, $G_2$ as follows.
\begin{itemize}
\item To build $G_1$, we apply the second part of Proposition \ref{propneed2} to $(\Theta_1,\,\rho,\,\ov\cM)$ in the role of $(\Theta,\,\rho,\,\cW)$, to get, in a computationally efficient way, a decomposition $\Theta_1=\sum_{i=1}^m\gamma_ig_{1,i}g_{1,i}^T$ with $g_{1,i}\in\ov \cM$ and $\gamma_i\geq0$, $\sum_i\gamma_i\leq 2\sqrt{2}\ln(4m^2n)\rho$. We set $G_1=[g_{1,1},...,g_{1,m}]$.
\item To build $G_2$, we subject $\Theta_2$ to the eigenvalue decomposition $\Theta_2=\Gamma\Diag\{\chi\}\Gamma^T$ and set $G_2=[g_{2,1},...,g_{2,m}]=Q\Gamma$.
\end{itemize}
Note that $\Theta_1+Q\Theta_2Q=\sum_i\gamma_ig_{1,i}g_{1,i}^T+\sum_i\chi_ig_{2,i}g_{2,i}^T$.

For the resulting polyhedral estimate $\wh w_{G,H}$ and for all $x_*\in\cX$, $s$-sparse $\nu_*$, and $\xi\in\Xi_\epsilon(G,H)$ if holds
\begin{equation}
\|\wh w_{G,H}(Ax_*+N\nu_*+\xi)-Bx_*\|\leq \ov f_H(\lambda,\mu,\tau,\Theta_1,\Theta_2,\rho)
\end{equation}
implying that the $\epsilon$-risk of the estimate is upper-bounded by $f_H(\lambda,\mu,\tau,\Theta_1,\Theta_2,\rho)$ (as $\xi\in\Xi_\epsilon(G,H)$
with probability $\geq1-\epsilon$).
\end{theorem}
Proof of the theorem follows that of Theorem \ref{the:newsparse} and is omitted.

\subsection{Putting things together}\label{together}
Finally, we can ``aggregate'' polyhedral estimates from Sections \ref{pgc1} and \ref{sectaltern} in the following construction (cf. \cite[Section 5.1.6]{PUP}):
\begin{quote}
Let us put
\[
\ov\varkappa(\epsilon)=\sigma\left(2\ln[(2n+2M+8m)/\epsilon]\right)^{1/2},
\] and let
$\wt H=[\wt h_1,...,\wt h_n]\in \bR^{m\times n}$, $\wt G=[\wt g_1,...,\wt g_{2m}]\in \bR^{m\times 2m}$, and $\ov H=[\ov h_1,...,\ov h_M]\in \bR^{m\times M}$, $\ov G=[\ov g_1,...,\ov g_{2m}]\in \bR^{m\times 2m}$ be the contrast matrices specified according to the synthesis recipes of Sections \ref{pgc1} and \ref{sectaltern}, respectively.
 We
define the aggregated estimate  $\wh w$ of $w_*$ as $\wh w(\omega)=B\wh x(\omega)$ where $\wh x(\omega)$ is the $x$-component of
\begin{align*}
(\widehat{x}(\omega),\widehat{\nu}(\omega))\in\Argmin_{x,\nu}
&\Big\{\|{\nu}\|_1: \;x\in\cX,\\
&\left.\begin{array}{ll}
|\wt h_k^T[{N\nu}+Ax-\omega]|\leq \ov\varkappa(\epsilon)\|\wt h_k\|_2,&k=1,...,n,\\
|\ov h_k^T[{N\nu}+Ax-\omega]|\leq \ov\varkappa(\epsilon)\|\ov h_k\|_2,&k=1,...,M,\\
|\wt g_i^T[{N\nu}+Ax-\omega]|_\infty\leq \ov\varkappa(\epsilon)\|\wt g_i\|_2,&i=1,...,2m,\\
|\ov g_i^T[{N\nu}+Ax-\omega]|_\infty\leq \ov\varkappa(\epsilon)\|\ov g_i\|_2,&i=1,...,2m,
\end{array}\right\}
\end{align*}
when the problem is feasible, and $\widehat{x}(\omega)=0$ otherwise.
\end{quote}
Let us denote $G=[\wt G,\ov G]\in \bR^{m\times 4m}$, let also $(\wt \lambda,\wt\mu,\wt\gamma)$ be a feasible solution to
the problem \rf{f_G2+} with $H=\wt H$, and let
$(\ov \lambda,\ov\mu,\ov\gamma)$ be a feasible solution to the problem \rf{f_G2++} with $H=\ov H$.
Let $f_{G,H}$ and $\bar f_{G,H}$ be specified in \rf{f_G2+} and \rf{f_G2++} respectively. From Propositions \ref{propro} and \ref{proproplus} it immediately follows  that for every $s$-sparse $\nu_*$ and every $x_*\in\cX$ the error bound
\beq\label{final}
\Risk_{\epsilon}[\widehat{w}(\cdot)|\cX,\cN]\leq \min[f_{G,\wt H}(\wt{\lambda},\wt{\mu},\wt\gamma),\bar f_{G,\ov H}(\ov\lambda,\ov\mu,\ov\gamma)]
\eeq
holds true.

Note that the resulting estimate can be efficiently optimized w.r.t. all parameters involved, except for $\ov H$, by specifying \begin{itemize}
\item $\wt H$ as (near) minimizer of $\rho[H]$ over $H\in \cH$  \rf{quick4},
\item  $\wt G$ as a result of the decomposition of the $(\Theta_1,\Theta_2)$-component of a (near-) optimal solution to the problem
(\ref{f_The22}), (\ref{f_The2f}) (see Theorem \ref{the:newsparse})
 associated with $\wt H$,
\item $\ov G$  as a result of the decomposition of the $(\Theta_1,\Theta_2)$-component of a (near-) optimal solution to the problem
(\ref{45old}) (see Theorem \ref{the:newsparseo}) associated with $\ov H$.
    \end{itemize}
\subsection{Numerical illustration}
\aic{}{In our ``proof of concept'' experiment we compare three estimates of $x_*$: 1) estimate $\wh x_{HG}$ with contrast matrix $[H,G]$ computed according to the recipe of Section \ref{pgc1}, 2) estimate $\wh x_{IG}$  with contrast $[\ov H,\ov G]=[I_m,\ov G]$ with $\ov G$ conceived utilizing
the synthesis routine of Section \ref{synt}, and 3) ``aggregated'' estimate $\wh x_{HIG}$ with combined contrast $[H,I_m,G,\ov G]$. }
\aic{In the above estimate}{We solve adopted versions of optimization problems in \rf{hko} and \rf{f_The22},  \rf{f_The2f} to compute contrasts $H$ and $G$ of the estimate $\wh x_{HG}$, and solve \rf{45old}, to build the contrast $\ov G$ of $\wh x_{IG}$.} For instance, when computing the contrast $\ov G$, we set $\ov \varkappa(\epsilon)= \sqrt{2}\sigma \mathrm{erfcinv}\left(\tfrac{\epsilon}{2n}\right)$ where $\mathrm{erfcinv}(\cdot)$ is the inverse complementary Gaussian error function; when processing problem (\ref{45old}) numerically,  $\Theta_1$ was set to 0; the resulting problem can be rewritten as
\begin{align}
\Opt=\min\limits_{\lambda,\mu,\gamma,\Theta}
&\Big\{\lambda+\aic{4\sum_{k=1}^{2n}\phi_\cS(\lambda)}{4\sum_{k=1}^\aic{n}{p} \mu_k+\sum_{j=1}^n \gamma_j}+
\aic{\rho_2}{}\Tr(\Theta):\,\lambda\geq0,\,\mu\geq0,\,\gamma\geq0,\,\Theta\succeq0,\nn
&\left.
\left[\begin{array}{c|c}\lambda \aic{I_n}{I_p}&\half \aic{B}{I_p}\cr\hline\half \aic{B^T}{I_p}&
\ov A^T \Theta \ov A+\aic{[P,\overline{A}]^T\Diag\{[\mu;\gamma]\}[P,\overline{A}]}{P^T\Diag\{\mu\}P+\overline{A}^T\Diag\{\gamma\}\overline{A}
}\end{array}\right]\succeq0
\right\}\label{13quick222}
\end{align}
where $\overline{A}=\aic{2}A/\rho_2$ with $\rho_2=2\ov \varkappa(\epsilon)+\varrho_{2,\ov H}$, the subsequent entries in $Pz$ being $z_1$, $[z_2-z_1]/h$, $[z_{i-2}-2z_{i-1}+z_i]/h^2$,
$3\leq i\leq \aic{n}{p}$, and  $h=2\pi/\aic{n}{p}$. \aic{Problem (\ref{13quick222}) was processed by the first order Composite Truncated Level method \cite{CTBL} in the same
fashion as problem (14) in \cite{CTBL}.}{ The corresponding risk bounds are evaluated by computing solutions to \rf{f_G2++}. Optimization problems involved are processed using {\tt Mosek} commercial solver \cite{mosek} via {\tt CVX} \cite{cvx2014}.} \footnote{MATLAB code for this  experiment is available at GitHub repository
 \url{https://github.com/ai1-fr/poly-robust}.}
\par
\begin{figure}[!htb]
\begin{tabular}{cc}
\includegraphics[width=0.45\textwidth]{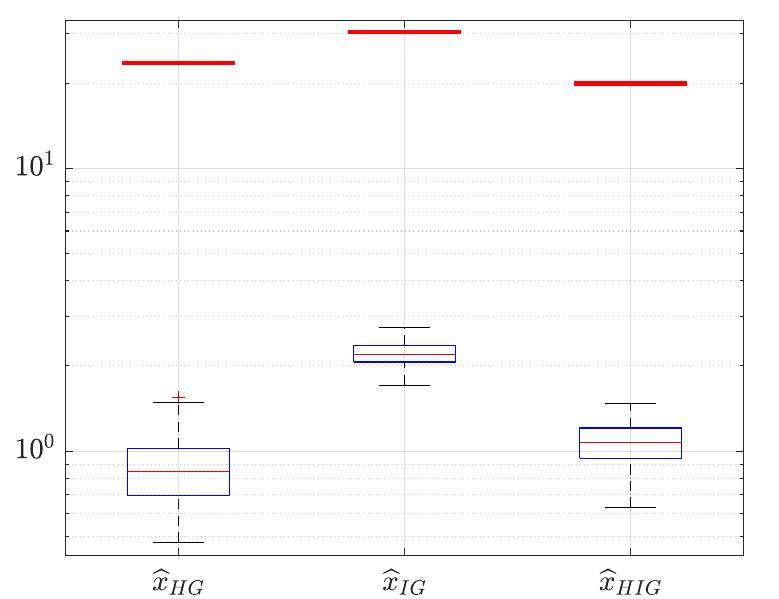}&
\includegraphics[width=0.45\textwidth]{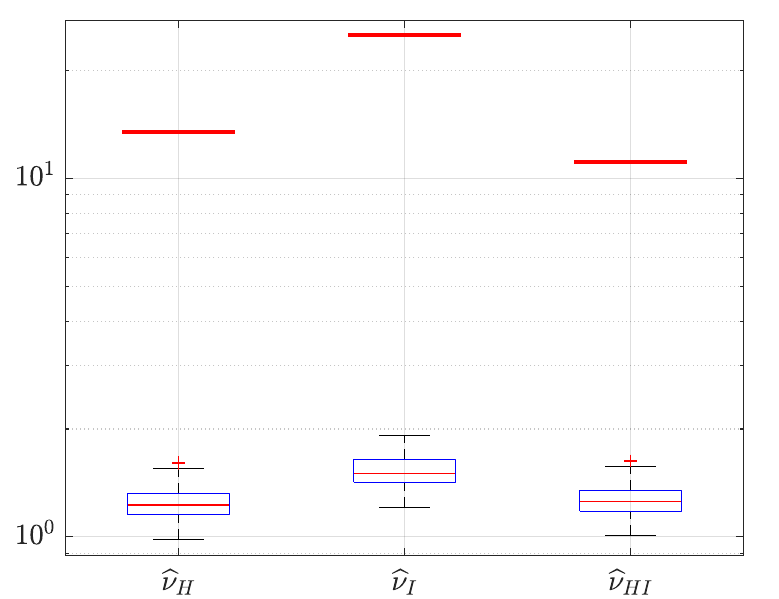}
\end{tabular}
\caption{\label{f:1} Left plot: distributions of $\|\cdot\|_2$-errors of recovery of $x_*$
and theoretical upper bounds on $\Risk_{0.05}$ (red horizontal bars); right plot: distributions of $\|\cdot\|_2$-errors
and theoretical upper bounds on $\Risk_{0.05}$  of recovery of $\nu_*$.}
\end{figure}
\begin{figure}[!htb]
\begin{tabular}{cc}
\includegraphics[width=0.45\textwidth]{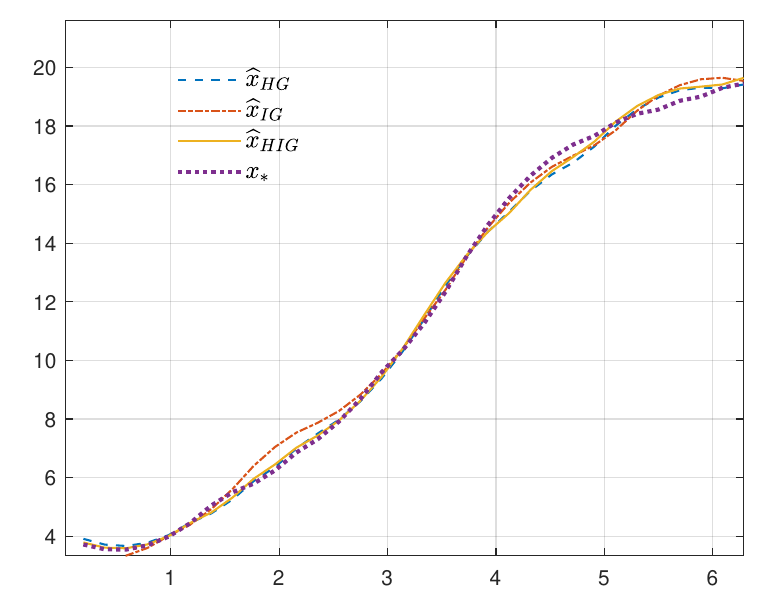}&
\includegraphics[width=0.45\textwidth]{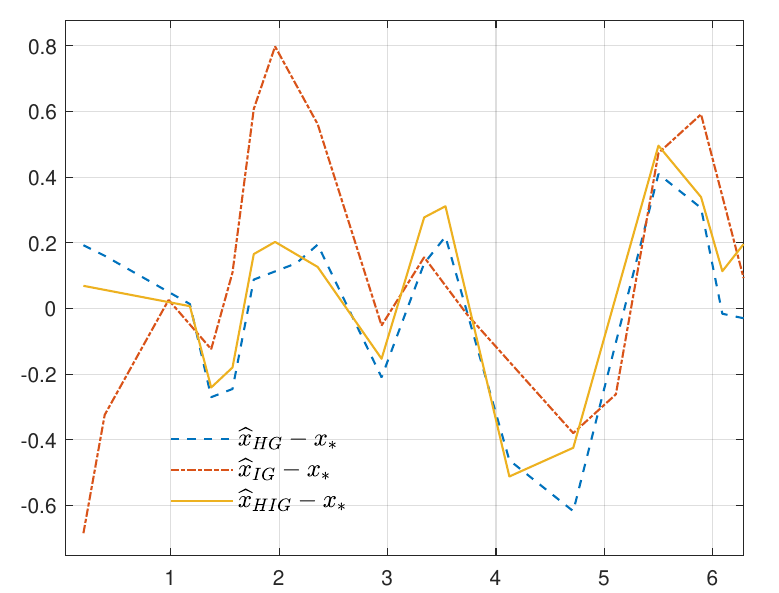}
\end{tabular}
\caption{\label{f:2} A typical signal/estimates realization and recovery errors.}
\end{figure}

In our illustration,
\begin{itemize}
\item $m=n=256$, $q=p=32$, $N=I_n$, $B=I_p$, $A$ is a $n\times p$ random matrix with Gaussian entries such that $A^TA=I_p$;
\item $\cX$ is the restriction on the $p$-point equidistant grid on the segment $\Delta=[0,2\pi]$ of functions $f$ satisfying
$|f(0)|\leq 4$, $|f'(0)|\leq1$, $|f''(t)|\leq4$, $t\in\Delta$;
\item the norm $\|\cdot\|$ quantifying the recovery error is the standard Euclidean norm on $\bR^p$;
\item $\xi\sim \cN(0,\sigma^2 I_m)$ with $\sigma=0.1$, $\epsilon=0.05$, and $s=8$.
\end{itemize}
Figure \ref{f:1} illustrates the results of the computation. In each experiment, we compute $n_S=100$ recoveries $\wh x_{HG}$, $\wh x_{IG}$, and $\wh x_{HIG}$ of randomly selected signals $x_*\in \cX$ with generated at random sparse nuisances $\nu_*$. The results are presented in the left plot in Figure \ref{f:1}. The right plot displays the boxplots of errors of recovery of the nuisance $\nu_*$ along with the upper risk bound $\Opt_2$ of \rf{Opt2max}. Figure \ref{f:2} illustrates a typical realization of the signal and the recovery errors; the  values of $\|\cdot\|_2$-recovery errors are $\|\wh x_{HG}-x_*\|_2=1.48...$, $\|\wh x_{IG}-x_*\|_2=2.02...$, and $\|\wh x_{HIG}-x_*\|_2=1.43...$, the corresponding $\|x_*\|_2=72.2...$.

\appendix
\section{Error  bound for $\ell_1$ recovery}\label{sectspars}
\paragraph{Condition $Q_\infty(s,\kappa)$}\begin{quote}
{\em Given an $m\times n$ sensing matrix $N$, positive integer $s\leq n$, and $\kappa\in(0,1/2)$, we say that $m\times p$ matrix $H$  satisfy condition $\bQ_\infty(s,\kappa)$ if}
\begin{equation}\label{ch1eq31}
\|w\|_{\infty}\leq \|H^TNw\|_\infty+\tfrac{\kappa}s \|w\|_1\quad\forall w\in\bR^n.
\end{equation}
\end{quote}

For $y\in\bR^n$, let $y^s$ stand for the vector obtained from $y$ by zeroing our all but the $s$ largest in magnitude entries.
\begin{proposition}\label{newl1theorem} Given $N$ and integer $s>0$, assume that $H$ satisfies the condition $\bQ_\infty(s,\kappa)$ with $\kappa<\half$.
Then for all $\nu,\widehat{\nu}\in\bR^n$ such that $\|\widehat{\nu}\|_1\leq\|{\nu}\|_1$ it holds:
\begin{equation}\label{napeq601}
\|\widehat{\nu}-{\nu}\|_q
\leq \tfrac{(2s)^{\frac{1}{q}}}{1-2\kappa}\left[\|H^TN[\widehat{\nu}-{\nu}]\|_\infty+\tfrac{\|\nu-{\nu}^s\|_1}{s}\right],\;1\leq q\leq \infty.
\end{equation}
\end{proposition}\noindent
{\bf Proof.} Let us denote
$
\rho=\|H^TN[\widehat{\nu}-{\nu}]\|_\infty,$
 and let
$z=\widehat{\nu}-{\nu}$.
\paragraph{1$^o$.} Let $I\subset \{1,...,n\}$ of cardinality $\leq s$ and let ${\ov I}$ be its complement in $\{1,...,n\}$. When denoting by $x_I$ the vector obtained from a vector $x$ by zeroing out the entries with indexes not belonging to $I$, we  have
\[
\|\widehat{\nu}_{\ov I}\|_1 = \|\widehat{\nu}\|_1-\|\widehat{\nu}_I\|_1\leq \|{\nu}\|_1-\|\widehat{\nu}_I\|_1
= \|\nu_I\|_1+\|\nu_{{\ov I}}\|_1-\|\widehat{\nu}_I\|_1\leq \|z_I\|_1+\|\nu_{{\ov I}}\|_1,
\]
and therefore
$$
\|z_{{\ov I}}\|_1\leq \|\widehat{\nu}_{{\ov I}}\|_1+\|\nu_{{\ov I}}\|_1\leq \|z_I\|_1+2\|\nu_{{\ov I}}\|_1.
$$
It follows that
\begin{equation}\label{napeq7}
 \|z\|_1=\|z_I\|_1+\|z_{{\ov I}}\|_1\leq 2\|z_I\|_1+2\|\nu_{{\ov I}}\|_1.
\end{equation}
 Besides this, by definition of $\rho$ we have
\begin{equation}\label{napeq10}
 \|H^TNz\|_\infty\leq \rho.
\end{equation}
\paragraph{2$^o$.}  Since $H$ satisfies $\bQ_\infty(s,\kappa)$, we have
\[
\|z\|_{s,1}\leq s\|H^TNz\|_\infty+\kappa \|z\|_1
\]
where $\|z\|_{s,1}$ is the $\ell_1$-norm of the $s$-dimensional vector composed of the $s$ largest in magnitude entries of $z$.
 By (\ref{napeq10}) it follows that $\|z\|_{s,1}\leq s\rho+\kappa \|z\|_1$ which combines with the evident inequality $\|z_I\|\leq \|z\|_{s,1}$ (recall that $\Card(I)=s$) and with (\ref{napeq7}) to imply that
\[\|z\|_1\leq 2\|z_I\|_1+2\|\nu_{{\ov I}}\|_1\leq 2s\rho+2\kappa\|z\|_1+2\|\nu_{{\ov I}}\|_1,
\]
hence (recall that $\kappa\leq \half$)
\be \|z\|_1\leq \frac{2s\rho+2\|\nu_{{\ov I}}\|_1}{1-2\kappa}.
\ee{napeq11}
On the other hand, since $H$ satisfies $\bQ_\infty(s,\kappa)$, we also have
\[\|z\|_\infty\leq \|H^TNz\|_\infty+\tfrac{\kappa}{ s}\|z\|_1,
\] which combines with (\ref{napeq11}) and (\ref{napeq10}) to imply that
\be
 \|z\|_\infty\leq \rho+\frac{\kappa}{ s} \frac{2s\rho+2\|\nu_{{\ov I}}\|_1}{1-2\kappa}=(1-2\kappa)^{-1}\left[\rho+\tfrac{\|\nu_{{\ov I}}\|_1}s\right]
\ee{napeq12}

We conclude that for all $1\leq q\leq \infty$,
\[\|z\|_p\leq\|z\|_\infty^{\frac{q-1}{q}}\|z\|_1^{\frac{1}{q}}\leq \tfrac{(2s)^{\frac{1}q}}{1-2\kappa}\left[\rho+\tfrac{\|\nu_{{\ov I}}\|_1}s\right].\eqno{\mbox{\qed}}
\]

\end{document}